\documentclass[11pt]{article} 

\usepackage[utf8]{inputenc}

\usepackage{geometry} 
\usepackage{graphicx}

\usepackage{booktabs}
\usepackage{array} 
\usepackage{paralist} 
\usepackage{verbatim} 
\usepackage{subfig}

\usepackage{fancyhdr} 
\usepackage{sectsty}
\allsectionsfont{\sffamily\mdseries\upshape} 

\usepackage[nottoc,notlof,notlot]{tocbibind} 
\usepackage[titles,subfigure]{tocloft}

\usepackage{titlesec}

\titleformat{\section}    
       {\normalfont\fontfamily{phv}\fontsize{16}{19}\bfseries}{\thesection}{1em}{}

\titleformat{\subsection}    
       {\normalfont\fontfamily{phv}\fontsize{12}{17}\bfseries}{\thesubsection}{1em}{}
\usepackage{threeparttable}
\usepackage{colortbl}

\usepackage{bm}
\usepackage{dsfont}
\usepackage{enumitem}
\usepackage{amsmath,amssymb}
\usepackage{hyperref}
 \usepackage[capitalise,noabbrev]{cleveref}

\crefformat{equation}{(#2#1#3)}

\newcommand{\customstretch}{1.3}
\usepackage{graphicx}
\usepackage{algorithm}
\usepackage{algpseudocode}
\usepackage{multirow}
\usepackage{tikz}

\newtheorem{theorem}{Theorem}

\newtheorem{proposition}[theorem]{Proposition}

\newtheorem{remark}{Remark}

\newtheorem{proof}{Proof}
\crefformat{equation}{(#2#1#3)}
\crefrangeformat{equation}{(#3#1#4) to~(#5#2#6)}
\crefmultiformat{equation}{(#2#1#3)}{ and~(#2#1#3)}{, (#2#1#3)}{ and~(#2#1#3)}
\newcommand*\circNum[1]{\tikz[baseline=(num.base)]{\node[shape=circle,inner sep=1pt,fill=orange,text=white,font=\bf] (num) {#1};}}

\newcommand{\N}{{\mathbb{N}}}
\newcommand{\R}{{\mathbb{R}}}
\newcommand{\Cn}{{\mathbb{C}}}
\DeclareMathOperator*{\colspan}{colspan}
\newcommand{\abs}[1]{\left\vert#1\right\vert}

\newcommand{\rT}[1]{#1^{\textsf{T}}}
\newcommand{\rTb}[1]{\rT{\left(#1\right)}}

\newcommand{\si}[1]{{#1^{+}}}

\newcommand{\grad}[1][]{{\nabla_{#1}}}

\newcommand{\dd}[2][]{{\frac{\mathrm{d}^{#1}}{\mathrm{d}#2}}}
\newcommand{\ddfx}{{\dd{\fx}}}

\newcommand{\ddt}{{\dd{t}}}

\newcommand{\Jtwo}[1]{{\mathbb{J}_{2#1}}}
\newcommand{\TJtwo}[1]{{\rT{\mathbb{J}}_{2#1}}}
\newcommand{\Jtn}{{\Jtwo{n}}}

\newcommand{\JtN}{{\Jtwo{N}}}
\newcommand{\TJtN}{{\TJtwo{N}}}
\newcommand{\Jtm}{{\Jtwo{m}}}

\newcommand{\I}[1]{{\fI_{#1}}}

\newcommand{\Z}[1]{{\fzero}_{#1}}

\newcommand{\tInit}{{t_{\mathrm{0}}}}
\newcommand{\tEnd}{{t_{\mathrm{end}}}}
\newcommand{\It}{I_t}
 
\newcommand{\fxInit}{{\fx_{\mathrm{0}}}}

\newcommand{\paramDomain}{\mathcal{P}}
\newcommand{\np}{{n_{\mathrm{\fmu}}}}
\newcommand{\POD}{\mathrm{POD}}
\newcommand{\DEIMalg}{\texttt{DEIM\_idx}}

\newcommand{\reduce}[1]{#1_{\mathrm{r}}}
\newcommand{\fxr}{\reduce{\fx}}
\newcommand{\fxrInit}{{\fx}_{\mathrm{r,0}}}
\newcommand{\ffr}{{\reduce{\ff}}}

\newcommand{\redSpace}{{\mathcal{V}}}

\newcommand{\solMnf}{\mathcal{M}}
\newcommand{\ns}{{n_{\mathrm{s}}}}
\newcommand{\xs}{\fx^{\mathrm{s}}}

\newcommand{\fXs}{\fX_{\mathrm{s}}}
\newcommand{\Xs}{{X_{\mathrm{s}}}}

\newcommand{\idxsDEIM}{\mathcal{I}_{\mathrm{DEIM}}}

\newcommand{\linSpace}{\mathbb{V}}
\newcommand{\symplForm}{\omega}

\newcommand{\Ham}{\mathcal{H}}
\newcommand{\Hamr}{\reduce{\Ham}}

\newcommand{\fb}{{\bm{b}}}
\newcommand{\fc}{{\bm{c}}}

\newcommand{\fe}{{\bm{e}}}
\newcommand{\ff}{{\bm{f}}}
\newcommand{\fg}{{\bm{g}}}

\newcommand{\fp}{{\bm{p}}}
\newcommand{\fq}{{\bm{q}}}
\newcommand{\fr}{{\bm{r}}}

\newcommand{\fx}{{\bm{x}}}
\newcommand{\fy}{{\bm{y}}}

\newcommand{\fA}{{\bm{A}}}
\newcommand{\fB}{{\bm{B}}}

\newcommand{\fD}{{\bm{D}}}

\newcommand{\fF}{{\bm{F}}}
\newcommand{\fG}{{\bm{G}}}
\newcommand{\fH}{{\bm{H}}}
\newcommand{\fI}{{\bm{I}}}

\newcommand{\fM}{{\bm{M}}}

\newcommand{\fP}{{\bm{P}}}

\newcommand{\fS}{{\bm{S}}}

\newcommand{\fU}{{\bm{U}}}
\newcommand{\fV}{{\bm{V}}}
\newcommand{\fW}{{\bm{W}}}
\newcommand{\fX}{{\bm{X}}}
\newcommand{\fY}{{\bm{Y}}}
\newcommand{\fZ}{{\bm{Z}}}

\newcommand{\fmu}{{\bm{\mu}}}

\newcommand{\fxi}{{\bm{\xi}}}

\newcommand{\frho}{{\bm{\rho}}}

\newcommand{\fPhi}{{\bm{\varPhi}}}
\newcommand{\fPsi}{{\bm{\Psi}}}

\newcommand{\fzero}{\ensuremath{\bm{0}}}

\newcommand{\calC}{\mathcal{C}}

\newcommand{\dict}[1]{\fD_{\text #1}}
\newcommand{\dictX}{\dict{X}}
\newcommand{\dictF}{\dict{F}}
\newcommand{\dictP}{\dict{P}}
\newcommand{\NX}{N_{\text{X}}}

\newcommand{\NP}{N_{\text{P}}}

\newcommand{\eig}{\texttt{eig}}
\newcommand{\fAX}{\fA_\textrm{X}}
\newcommand{\xInitXs}{\fx_\textrm{0,X,s}}
\newcommand{\xInitX}{\fx_\textrm{0,X}}
\newcommand{\fGYs}{\fG_\textrm{Y,s}}
\newcommand{\fGXs}{\fG_\textrm{X,s}}

\newcommand{\fAXs}{\fA_\textrm{X,s}}

\newcommand{\fHX}{\fH_\textrm{X}}
\newcommand{\fHXJr}{\fH_{\textrm{X},\JtN, \textrm{r}}}
\newcommand{\fHXJl}{\fH_{\textrm{X},\JtN, \textrm{l}}}
\newcommand{\fHXJJ}{\fH_{\textrm{X},\JtN, \JtN}}
\newcommand{\fHXJrs}{\fH_{\textrm{X},\JtN, \textrm{r,s}}}
\newcommand{\fHXJls}{\fH_{\textrm{X},\JtN, \textrm{l,s}}}
\newcommand{\fHXJJs}{\fH_{\textrm{X},\JtN, \JtN, \textrm{s}}}

\newcommand{\fHXs}{\fH_\textrm{X,s}}

\newcommand{\fGX}{\fG_\textrm{X}}
\newcommand{\fGXF}{\fG_\textrm{X,F}}
\newcommand{\fGXFJ}{\fG_{\textrm{X,F},\JtN}}
\newcommand{\fGXFJs}{\fG_{\textrm{X,F},\JtN, \textrm{s}}}
\newcommand{\fGXJ}{\fG_{\textrm{X},\JtN}}
\newcommand{\fGXJs}{\fG_{\textrm{X},\JtN,\textrm{s}}}
\newcommand{\fGXFs}{\fG_\textrm{X,F,s}}
\newcommand{\fGF}{\fG_\textrm{F}}
\newcommand{\fGFs}{\fG_\textrm{F,s}}
\newcommand{\fUP}{\fU_{\fP}}
\newcommand{\fUPPo}{\fU_{\hat\fP, \fP_o}}
\newcommand{\fFhatP}{\fF_{\hat \fP}}
\newcommand{\fFhatPs}{\fF_{\hat \fP, \textrm{s}}}
\newcommand{\fGVU}{\fG_{\fV,\fU}}

\title{Dictionary-based Online-adaptive Structure-preserving Model Order Reduction for Parametric Hamiltonian Systems}
\author{R. Herkert\footnote{robin.herkert@ians.uni-stuttgart.de, Institute of Applied Analysis and Numerical Simulation, University of Stuttgart,Pfaffenwaldring 57, Stuttgart, 70569, Germany}, P. Buchfink\footnote{patrick.buchfink@ians.uni-stuttgart.de, Institute of Applied Analysis and Numerical Simulation, University of Stuttgart,Pfaffenwaldring 57, Stuttgart, 70569, Germany}, B. Haasdonk\footnote{bernard.haasdonk@ians.uni-stuttgart.de, Institute of Applied Analysis and Numerical Simulation, University of Stuttgart,Pfaffenwaldring 57, Stuttgart, 70569, Germany}}

\date{}
\begin{document}
\maketitle
\begin{abstract}
Classical model order reduction (MOR) for parametric problems may become computationally inefficient due to large sizes of the required projection bases, especially for problems with slowly decaying Kolmogorov $n$-widths. Additionally, Hamiltonian structure of dynamical systems may be available and should be preserved during the reduction. In the current presentation, we address these two aspects by proposing a corresponding dictionary-based, online-adaptive MOR approach. The method requires dictionaries for the state-variable, non-linearities and discrete empirical interpolation (DEIM) points. During the online simulation, local basis extensions/simplifications are performed in an online-efficient way, i.e.\ the runtime complexity of basis modifications and online simulation of the reduced models do not depend on the full state dimension. Experiments on a linear wave equation and a non-linear Sine-Gordon example demonstrate the efficiency of the approach.
\end{abstract}

approximation, energy preservation


\section{Introduction}
The success of classical Model Order Reduction (MOR) methods relies on the assumption that the set of solutions can be approximated well in a low-dimensional subspace.
However, many problems do not allow such a low-dimensional approximation,
which often leads to reduced bases that are too large
such that no sufficient speed-up can be realized compared to the high-order simulation.
These problems can be characterized by slowly decaying Kolmogorov $n$-widths \cite{Pinkus1985}.
The development of efficient MOR for such problems is still a strongly investigated topic.

One class of MOR techniques for problems with slowly decaying Kolmogorov $n$-widths are online-adaptive approaches
which are characterized by adapting the reduced basis during the online phase.
This online-adaption may be a low-rank update of the reduced basis \cite{Peherstorfer2020,Zimmermann2018}.
Another approach is local MOR which precomputes multiple reduced bases and switches to the appropriate basis depending on different notions of locality \cite{Amsallem2016,Amsallem2012,Dihlmann2011,Drohmann2011,Eftang2010}.
Furthermore, dynamical low-rank techniques adapt the basis by including additional equations for the evolution of the reduced basis on a manifold of low-rank matrices \cite{Koch2007}.

In the scope of this work,
we focus on dictionary-based online-adaptive MOR.
The idea is to build a dictionary from which a reduced basis is computed in the online phase.
In a suitable offline stage, the computational intensive operations, such as the computation of the corresponding reduced operators based on the dictionary, are precomputed and stored. 
This makes dictionary-based MOR online-efficient, i.e.\ the computational costs during the online phase do not depend on the state dimension of the full model.
Meanwhile, various extensions have been proposed for dictionary-based MOR:
A residual-based selection scheme for stationary elliptic problems has been introduced in \cite{kaul13},
while unsteady linear transport problems have been treated in \cite{dihl12}.
An approach for steady problems with a local anisotropic parameter-space distance measure as selection criterion has been given in \cite{stam13}.
For non-linear hyperbolic equations,
an $L^1$-norm minimization of the residual has been introduced \cite{Abgrall2016,Abgrall2018} together with an error estimator.
In \cite{Balabanov2021},
dictionary-based MOR is enhanced with methods from randomized linear algebra.
To this end,
a sketched sparse minimal residual approximation is introduced,
which reduces the computational cost and storage while showing increased numerical stability.

Another type of online-adaptivity in the context of MOR worthwhile mentioning adapts the approximation of the non-linearity in the scope of hyper-reduction during the online phase, e.g.\ in \cite{Peherstorfer2014,Peherstorfer2015}.

Furthermore, models may possess the form of Hamiltonian systems.
A Hamiltonian system has a certain structure
which ensures conservation of energy and, under mild assumptions, stability properties.
Classical MOR techniques like the Proper Orthogonal Decomposition (POD), see e.g.\ \cite{Volkwein2013}, fail to preserve this Hamiltonian structure, which, in general, violates the conservation of energy and may yield unstable reduced models.
To this end, structure-preserving MOR for Hamiltonian systems (also symplectic MOR) has been introduced \cite{MaboudiAfkham2017,Peng2016}.
It relies on the generation of a symplectic basis and a projection with the so-called symplectic inverse of the basis.
An adaption of the Discrete Empirical Interpolation Method (DEIM) \cite{Chaturantabut2010}, the symplectic DEIM (SDEIM), has been introduced in \cite{Peng2016} for the hyper-reduction of non-linear Hamiltonian systems.

In the present work,
we merge ideas from structure-preserving MOR for Hamiltonian systems with dictionary-based MOR and present online-efficient methods.
Our three key contributions are:
\begin{enumerate}
  \item We introduce dictionary-based structure-preserving MOR for Hamiltonian systems,
  \item to this end, we introduce a new symplectic basis generation technique, where the basis is computed from a dictionary, the dictionary-based complex SVD (DB-cSVD),
  \item we treat non-linearities online-efficiently with a new dictionary-based SDEIM (DB-SDEIM). 
\end{enumerate}

There is recent work which also considers online-adaptive structure-preserving MOR for Hamiltonian systems \cite{Pagliantini2021}.
Instead of dictionary-based techniques, it uses ideas from dynamical low-rank methods.
In contrast to our approach however, the computational costs in the online-phase of these methods depend linearly on the state dimension of the unreduced model.

This work is structured as follows:
In \Cref{sec:essentials} essential background on MOR and Hamiltonian systems will be given.
In \Cref{sec:main}, the new dictionary-based, structure-preserving methods will be presented.
\Cref{sec:numerics} will focus on numerical experiments and comparisons with standard structure-preserving methods.
The work is concluded in \Cref{sec:conclusion}.

\section{Essentials}\label{sec:essentials}
In this section,
we briefly introduce the reader to the essentials required for our study.
To this end, we discuss and introduce notation for classical MOR, 
Hamiltonian systems and symplectic MOR.

\subsection{Classical projection-based MOR}
\label{Subsec:classical_mor}

In the scope of this work, we aim to derive an efficient surrogate for a
high-dimensional, parametric dynamical system with MOR techniques.
In the following, we introduce this high-dimensional problem,
the projection-based reduced-order model,
and hyper-reduction of the reduced-order model for an efficient evaluation.
For a more detailed introduction to MOR, we refer to the monographs \cite{Benner2017,Benner2021}.

For a given \emph{parameter domain} $\paramDomain \subset \R^\np$,
\emph{time interval}
\footnote{Note that even the end time may be parametric ($\tEnd(\fmu)$) leading to parameter-dependent time intervals $\It(\fmu)$, as we use in the experiments. For notational simplicity however, we use a fixed time interval in the methods' description.}
$\It := [\tInit, \tEnd] \subset \R$,
\emph{right-hand side} (RHS) $\ff: \R^N \times \It \times \paramDomain \to \R^N$,
and \emph{initial value} $\fxInit: \paramDomain \to \R^N$,
the \emph{full-order model (FOM)} reads:
For a fixed (but arbitrary) \emph{parameter vector} $\fmu \in \paramDomain$
find the corresponding \emph{solution} $\fx(\cdot; \fmu) \in \calC^1(\It, \R^N)$
with
\begin{equation}\label{eq:FOM}
\begin{aligned}
  \ddt \fx(t; \fmu) &= \ff(\fx(t; \fmu), t; \fmu) \qquad \text{for all } t \in \It,\\
  \fx(\tInit; \fmu) &= \fxInit(\fmu).
\end{aligned}
\end{equation}
For well-posedness of the FOM, we assume that the RHS $\ff$
is continuous in $t$ and Lipschitz-continuous in $\fx$ and that $\tEnd$ is small enough
such that a unique solution of the FOM exists.

The set of all solutions
\begin{align*}
 \solMnf := \left\{
   \fx(t; \fmu) \,\vert  \, (t, \fmu) \in \It \times \paramDomain
 \right\} \subset \R^N
\end{align*}
is the so-called \emph{solution manifold}.
Classical projection-based MOR proceeds in two steps to approximate the solution manifold:
Firstly, the solution is approximated in a low-dimensional subspace 
$\redSpace = \colspan\left({\fV}\right) \subset \R^N$ with
$\dim\left(\redSpace\right) = n$, $n \ll N$ and
\begin{align*}
  \fx(t; \fmu) \approx \fV \fxr(t; \fmu)
\qquad \text{for all } (t, \fmu) \in \It \times \paramDomain
\end{align*}
for some \emph{reduced coefficients} $\fxr(t; \fmu) \in \R^{n}$
via the reduced-order basis (ROB) matrix $\fV \in \R^{N \times n}$.
Secondly, the residual $\fr(t; \fmu)$ is required to vanish on
$\colspan\left( \fW \right)$
for the so-called \emph{projection matrix} $\fW \in \R^{N \times n}$,
\begin{align*}
  &\fr(t; \fmu) = \fV \ddt\fxr(t; \fmu) - \ff(\fV \fxr(t; \fmu), t; \fmu) \in \R^{N},&
  &\rT\fW \fr(t; \fmu) \stackrel{!}{=} \Z{n \times 1}.
\end{align*}
If the basis matrices are biorthogonal, $\rT\fW \fV = \I{n}$,
the resulting reduced-order model (ROM) reads:
For a fixed (but arbitrary) parameter vector $\fmu \in \paramDomain$,
find the corresponding \emph{reduced solution} $\fxr(\cdot; \fmu) \in \calC^1(\It, \R^N)$
with
\begin{equation}\label{eq:ROM}
\begin{aligned}
  \ddt \fxr(t; \fmu) &= \ffr(\fxr(t; \fmu), t; \fmu) \qquad \text{for all } t \in \It,\\
  \fxr(\tInit; \fmu) &= \fxrInit(\fmu)
\end{aligned}
\end{equation}
with the \emph{reduced right-hand side}
$\ffr(\fxr, t; \fmu) := \rT\fW \ff(\fV \fxr, t; \fmu)$
and the \emph{reduced initial value}
$\fxrInit(\fmu) := \rT\fW \fxInit(\fmu)$.

One class of techniques to determine a suitable reduced space
are \emph{snapshot-based basis generation techniques}.
These techniques consider a finite subset
$\Xs \subset \solMnf$, $\abs{\Xs} = n_{X_s} < \infty$, of the solution manifold.
The elements $\xs_i \in \Xs$, $1 \leq i \leq n_{X_s}$, are typically referred to
as \emph{snapshots}.
For convenience, we define the
\emph{snapshot matrix} $\fXs := (\xs_i)_{i=1}^{n_{X_s}} \in \R^{N \times n_{X_s}}$ as the matrix
which stacks the snapshots in its columns.
The most popular snapshot-based basis generation technique is the
POD.
It generates a reduced basis from the $m$ most important left-singular vectors of the snapshot matrix $\fXs$.
We denote this technique with $\fV = \POD(\fZ_s, m)$ where $\fZ_s \in \Cn^{N \times n_{Z_s}}$ is a, possibly complex, snapshot matrix.

If the RHS $\ff$ is linear and parameter-separable,
the ROM can be solved online-efficiently with a so-called offline--online decomposition
(see e.g.\ \cite{Haasdonk2011b}).
However, if the RHS is linear but not parameter-separable or non-linear,
another approximation is required for an efficient evaluation of the ROM.
These techniques are typically referred to as hyper-reduction.
We use the (Discrete) Empirical Interpolation Method ((D)EIM)
\cite{Barrault2004,Chaturantabut2010} for this purpose.  
The assumption of DEIM is
(i) that single components
$f_i := (\ff)_i$ for $1 \leq i \leq N$ of the RHS $\ff$
only depend on a few components of the state $\fx$
and (ii) that single components of the RHS can be evaluated efficiently.
The idea is to select a few important components $f_i$ which are computed
for a given index set of \emph{DEIM (interpolation) indices} 
$\idxsDEIM \subset \{1, \dots, N\}$
and the other components $f_i$ with $i \in \{1, \dots, N\} \setminus \idxsDEIM$
are interpolated.

An algorithm to determine a vector of DEIM indices and a corresponding projection matrix
is is reproduced from \cite{Chaturantabut2010} in Algorithm \Cref{DEIMalg}.
\begin{algorithm}
\caption{$\DEIMalg$-algorithm}\label{DEIMalg}
\textbf{Input}: Basis matrix $ \fU \in \R^{n\times m}$\\
\textbf{Output}: Index vector $\frho = (\rho_1,...,\rho_m\rT) \in \mathbb{N}^m$, proj.\ matrix $\fP \in \{0,1\}^{n\times m}$
\begin{algorithmic}[1]
\State$\fP= \Z{n\times m}$, $\frho = \Z{m\times 1}$
\State$[M,\rho_1]$ = max(abs($\fU(:,1)$)) \Comment{maximum and the corresponding index}
\State$\fP(\rho_1,1)$ = 1
\For {$\ell = 2:m$}
\State $\fc$ = $(\fP(:,1:\ell-1\rT)\fU(:,1:\ell-1))^{-1}(\fP(:,1:\ell-1\rT)\fU(:,\ell))$
\State$[M, \rho_\ell]$ = max(abs($\fU(:,\ell) - \fU(:,1:\ell-1)\fc$)) \Comment{max.\ and corresp.\ idx}
\State$\fP(\rho_\ell,\ell) = 1$
\EndFor
\end{algorithmic}
 \end{algorithm}

\subsection{Hamiltonian Systems and Symplectic MOR}
A special class of dynamical systems are Hamiltonian systems,
which are used to model energy-conserving processes.
In the following, we briefly discuss symplectic vector spaces,
Hamiltonian systems, and structure-preserving MOR for
parametric high-dimensional Hamiltonian systems.

A \emph{symplectic (vector) space} is a tuple $(\linSpace, \symplForm)$ of a
\emph{phase space} $\linSpace$
and a skew-symmetric and non-degenerate bilinear form
$\symplForm: \linSpace \times \linSpace \to \R$ called the \emph{symplectic form}.
It can be shown \cite{daSilva2008}
(i) that the phase space is necessarily even-dimensional which is why
we restrict to $\linSpace \equiv \R^{2N}$
and (ii) that there exist canonical coordinates such that the symplectic form
can be represented by the \emph{canonical Poisson tensor}
\begin{align*}
  \JtN := \begin{bmatrix}
    \Z{N} & \I{N}\\
    -\I{N} & \Z{N}
  \end{bmatrix} \in \R^{2N \times 2N}
\end{align*}
with the identity and zero matrix $\I{N}, \Z{N} \in \R^{N \times N}$.
Thus, we only consider the symplectic vector space $(\R^{2N}, \JtN)$
where we denote the symplectic form by its coordinate matrix $\JtN$.

A matrix $\fB \in \R^{2N \times 2m}$ defining a linear mapping between two symplectic spaces $(\R^{2m}, \Jtm)$ and $(\R^{2N}, \JtN)$ with $m \leq N$
is called a \emph{symplectic matrix} if
\begin{align*}
  \rT\fB \JtN \fB = \Jtm.
\end{align*}
Moreover, a differentiable function $\fg: (\R^{2m}, \Jtm) \to (\R^{2N}, \JtN)$
with $m \leq N$
is called a \emph{symplectic mapping}, if
the \emph{Jacobian} $\ddfx \fg(\fx) \in \R^{2N \times 2m}$ is
a symplectic matrix for all $\fx \in \R^{2N}$.

A \emph{Hamiltonian system (on a symplectic vector space in canonical coordinates)} is
a triple $(\R^{2N}, \JtN, \Ham)$ consisting of a symplectic vector space
$(\R^{2N},\JtN)$ with a canonical Poisson matrix $\JtN$ and a \emph{Hamiltonian
(function)} $\Ham(\cdot; \fmu) \in \calC^1(\R^{2N})$, which we additionally assume to be dependent on a parameter $\fmu \in \paramDomain$. 
Moreover, we assume to be given a parameter-dependent initial value $\fxInit: \paramDomain \to \R^{2N}$.
The parametric Hamiltonian system reads: for a fixed (but arbitrary) parameter vector
$\fmu \in \paramDomain$, find the solution $\fx(\cdot; \fmu) \in \calC^1(\It, \R^{2N})$
with
\begin{equation}\label{eq:ham_sys}
\begin{aligned}
  \ddt \fx(t; \fmu)
&= \JtN \grad[\fx] \Ham(\fx(t; \fmu); \fmu) \qquad \text{for all } t \in \It,\\
  \fx(\tInit; \fmu) &= \fxInit(\fmu).
\end{aligned}
\end{equation}
The most important property of a Hamiltonian system is that the solution conserves the Hamiltonian over time,
$\ddt \Ham(\fx(t; \fmu); \fmu) = 0$ for all $t \in \It$.
In some cases it is convenient to split the solution
$\fx(t; \fmu) = [\fq(t; \fmu); \fp(t; \fmu)]$ in separate coordinates
$\fq(t; \fmu),\; \fp(t; \fmu) \in \R^N$.

Symplectic MOR \cite{Peng2016,MaboudiAfkham2017} is a projection-based MOR technique
(see \Cref{Subsec:classical_mor}) to reduce parametric high-dimensional Hamiltonian systems.
The resulting ROM is a low-dimensional
Hamiltonian system $(\R^{2n}, \Jtn, \Hamr)$ which preserves
the so-called \emph{reduced Hamiltonian}
$\Hamr(\fxr) := \Ham(\fV \fxr)$.
This is ensured if (i) the ROB matrix $\fV$ is a symplectic matrix
and (ii) the projection matrix $\fW$ is set to be the transpose of the so-called \emph{symplectic inverse}
$\si{\fV}$ of the ROB matrix $\fV$, i.e.
\begin{align*}
  \rT\fW := \si{\fV} := \Jtn \rT\fV \TJtN.
\end{align*}

\section{Dictionary-based symplectic MOR}
\label{sec:main}
Standard, projection-based MOR requires that the solution manifold can be well approximated in a low-dimensional subspace. In technical terms, this means that the solution manifold has so-called rapidly decaying Kolmogorov $n$-widths \cite{Pinkus1985}. For many transport problems, however, this is not the case and thus the resulting global basis matrices $\fV, \fW$ require large $n$, hence the ROMs lose their efficiency for such problems. In these situations, dictionary-based approaches may still yield online-efficient reduced models.
The idea of a dictionary-based approach is that the basis is computed in the online-phase during the time stepping of the reduced simulation and is thus parameter- and time-dependent. Therefore, it can be expected that the basis may be much smaller and a higher speed-up may be realized compared to a standard approach, depending on how efficient the basis changes and basis computations can be performed \cite{dihl12,kaul13,stam13}.
To compute a basis, a dictionary of $\NX \in \N$ state snapshots 
\[
  \dictX := \{\xs_1, ..., \xs_{N_\text{X}}\} \subset \solMnf
\]
is constructed in the offline-phase and during the online-phase snapshots are selected from the dictionary, where \emph{dictionary} denotes a finite subset of $\R^{2N}.$
Then, the basis is computed from the selected snapshots.

We assume our dictionary to be labeled with a parameter-time-label, which allows us to use these labels for the selection rule. By labelling we mean, that we have a bijection $l:\dictX\to L\subset \R^{n_L}$, where $L$ is the set of labels. Typical choices for labels are time labels ($L \subset \It$), parameter labels ($L \subset \paramDomain$) or combinations of this, e.g.\ parameter-time-labels ($L = L_1 \times L_2$, $L_1 \subset \It$, $L_2 \subset \paramDomain$) which will be used in the following.

\subsection{Workflow}
\begin{figure}[H]
 \centering
 \includegraphics[width=0.8\linewidth]{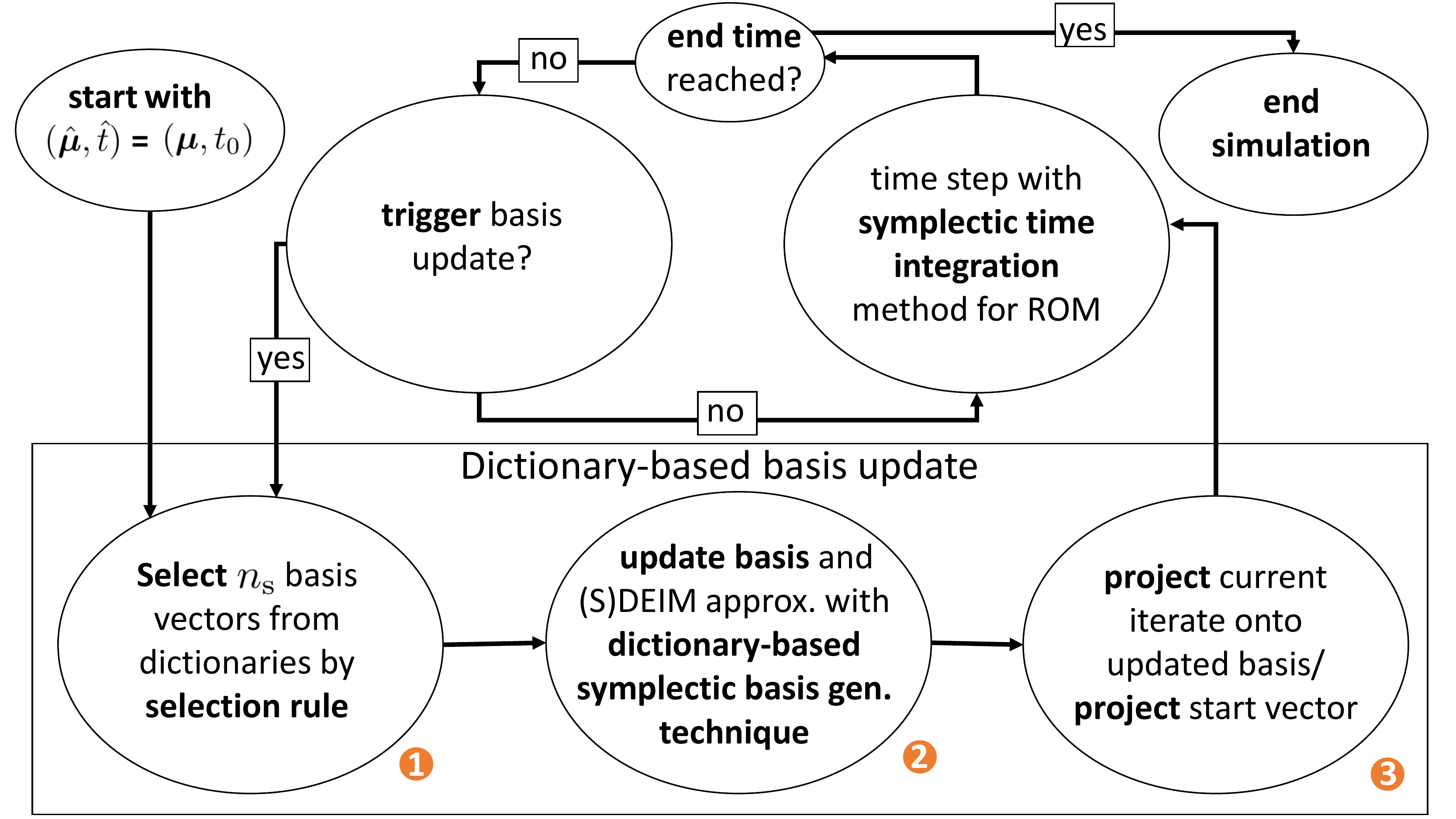}
 \caption{Workflow sketch for the online-phase of symplectic dictionary-based MOR}
 \label{fig:Workflow}
\end{figure}

In this section, the workflow of our dictionary-based symplectic MOR-algorithm is explained. In \Cref{fig:Workflow}, a workflow sketch is presented that describes the interplay of the methods described in the subsequent subsections. The algorithm starts with the selection of snapshots from the dictionary for a given parameter $\fmu$ and start time $\tInit$. Next, the basis and SDEIM approximation are computed using the dictionary-based methods. Then, the start vector is projected. After that, the time stepping for the reduced system is advanced until a basis update is triggered. We use a fixed trigger which means that we introduce a window size $m_s$ that specifies after how many time steps the basis should be updated. As soon as a basis update has been performed, the current state is projected from the old to the new reduced space and the time-stepping is continued. For the time integration it is recommended to use so-called symplectic integrators, which preserve the symplectic structure after time-discretization. As a symplectic integrator we use the implicit midpoint rule. These steps are repeated until we reach our final simulation time $\tEnd.$ The projections with the basis matrices are performed in an implicit, online-efficient manner. How this implicit formulations for different methods are obtained is explained in the Subsections \ref{secDPOD} to \ref{secDSMOR}. The selection procedure is described in \Cref{secOnlineSelection}. 

\subsection{Online-selection}
\label{secOnlineSelection} 
Firstly, we want to discuss, how snapshots can be selected from the dictionary with a selection rule (step~\circNum{1} of the dictionary-based basis update from \Cref{fig:Workflow}).
In \cite{stam13} for stationary problems, $\ns$ snapshots are selected during the online-phase, whose parameter-labels are closest to the one of the currently queried parameter. In the present paper, we extend this idea to time-dependent problems: An additional time-label is treated as an additional entry of the parameter vector to form our label. Furthermore, the window size $m_s\in\N$ is introduced which specifies how many time steps are to be performed with the same basis

During the reduced simulation the $\ns$ dictionary snapshots are selected, 
whose parameter-time-labels are closest to the set of parameter-time-labels of the current parameter and the current and next $m_s$ time steps:
$$ I_s=\{i_1,.., i_{n_s}\} = \textrm{argmin} (\{d_1,...,d_{N_{\textrm{X}}}\}, \ns)$$ with $d_{i} = \min(\{d_{i,0}, ...,d_{i,m_s}\}),\ i = 1,...,N_\textrm{X}$ and $$d_{i,\ell} = d_{\mathcal{P}T} (l(\xs_i),(\fmu,t+\ell\Delta t)),\ \ell = 0,...,m_s, $$
for a metric $d_{\mathcal{P}T}:L\times L \to \R^+_0.$ The procedure is presented in \Cref{Selalg}.

\begin{algorithm}
  \caption{\texttt{DB\_indices\_selection}}\label{Selalg}
  \textbf{Input}: Window size $m_s$, current time $t$, time step width $\Delta t$, number of snapshots to select $\ns$, parameter vector $\fmu$, metric $d_{\mathcal{P}T}$, label function $l$\\
  \textbf{Output}: Set of indices $I_s$ from which the next basis is computed 
  \begin{algorithmic}[1]
    \For{$i = 1,.., N_{\textrm{X}}$}
    \For{$\ell = 0,...,m_s$}
    \State$d_{i,\ell} = d_{\mathcal{P}T} (l(\xs_i),(\fmu,t+\ell\Delta t))$ \Comment{compute distances between labels}
    \EndFor
    \State $d_{i} = \min(\{d_{i,0}, ...,d_{i,m_s}\})$ \Comment{Comp.\ min.\ dist.\ for dictionary elements}
    \EndFor
    \State $ I_s=\{i_1,.., i_{n_s}\} = \textrm{argmin} (\{d_1,...,d_{N_{\textrm{X}}}\}, \ns)$ \Comment{sel.\ idxs.\ with smallest dists}
  \end{algorithmic}
\end{algorithm}
As a metric, a standard 2-norm $d_{\mathcal{P}T}((\fmu,t),(\hat \fmu,\hat t)) = \|(\fmu,t) - (\hat \fmu,\hat t)\|$ could be used. We use a weighted 2-norm with scaled time-distance 
$$d_{\mathcal{P}T}((\fmu,t),(\hat\fmu,\hat t)) =  \sqrt{\|\fmu - \hat \fmu\|_2^2+c\vert t - \hat t\vert^2},$$ 
with $c = \frac{d_{L, \textrm{max}}^2 }{\Delta t^2}$ because usually the time-step width will be much smaller than the distance between the parameters from which the dictionary is sampled. The quantity $d_{L, \textrm{max}}$ is defined as 
$$d_{L, \textrm{max}} = \max\limits_{\fmu \in L}\min\limits_{\hat\fmu \in L}\|\fmu-\hat \fmu \|_2.$$
We formally select snapshots from the state dictionary by defining the global snapshot matrix $\fX = [\xs_1, \dots, \xs_{\NX}] \in \R^{2N \times \NX}$ and a selection matrix $\fP_s \in \{0,1\}^{\NX\times \ns},$ where $\ns$ denotes the number of selected snapshots, and forming the product $\fX \fP_s$. The matrix $\fP_s$ is not actually computed, but the formal multiplications with it or its transposed are obtained by corresponding column and row selections.

\subsection{Dictionary-based POD}
\label{secDPOD}
As a possibility for a dictionary-based basis computation (step~\circNum{2} and ~\circNum{3} of the dictionary-based basis update from \Cref{fig:Workflow}), we discuss the dictionary-based POD (DB-POD), which is compared to the standard POD in Table \ref{tablestandvsdictbasedPOD}. The  corresponding mathematical operations of these two procedures are aligned with each other. For notational simplicity we restrict to the case of a linear and non-parametric system with time-independent RHS, i.e.\ $ \ff(\fx(t; \fmu), t; \fmu) = \fA\fx$, $\fxInit(\fmu) = \fxInit$. Additional separable parameter dependence of $\fA$ and $\fxInit$ can be allowed by performing the operations to the components of $\fA$ and $\fxInit$ and then performing an additional online linear combination step. 
The basic idea of dictionary-based techniques is that during an offline-phase high-dimensional matrix products are precomputed once and then are used during the online-phase to assemble the reduced system efficiently. The algorithm is obtained by collecting the matrix products which are necessary for the computation of the POD-basis and using the fact that the POD basis $\fV := \POD(\fX \fP_s, m)$ can be expressed by a multiplication of the snapshot matrix with the eigenvectors of $\rT\fP_s\rT\fX \fX \fP_s\in \R^{\ns\times\ns}$ and scaling with the inverses of the square roots of its eigenvalues $\fV :=\fX \fP_s \fPhi\fS $ with $\fS = $diag$\left(\left(\frac{1}{\sqrt{\lambda_1}}, ... ,\frac{1}{\sqrt{\lambda_{m}}}\right)\right) \in \R^{n_s \times m}$ with $\lambda_1\geq... \geq\lambda_{m}>0$ the dominant, non-zero eigenvalues of $\rT\fP_s\rT\fX \fX \fP_s.$ 
Like this, the projections with the POD basis can be computed as a multiplication of low-dimensional matrices, which are selected as sub-matrices of precomputed quantities and low-dimensional vectors, computed from a low-dimensional eigenvalue problem.  
\begin{table}[h]
\caption{Comparison of dictionary-based POD vs.\ standard POD for a linear non-parametric system} 
\label{tablestandvsdictbasedPOD}
\centering
{\renewcommand{\arraystretch}{\customstretch}%
\small
\begin{tabular}{|c|c|}
\hline
\bf Dictionary-based POD: & \bf Standard POD(of $\fX\fP_s$):\\
\hline
\bf Offline-DB-POD: &\\
\hline
Glue all $\fx(t_j, \fmu_i)\in \dictX$ together as $\fX$, &\\
precompute $\fG_\textrm{X} := \rT\fX \fX$, $\fA_\textrm{X} := \rT\fX\fA \fX$,&\\ 
$\fx_\textrm{0,X} := \rT\fX \fxInit$ & \\
\hline
\bf Online-DB-POD: & \\
\hline
Select sub-matrix $\fG_\textrm{X,s} := \rT\fP_s\fG_\textrm{X} \fP_s$ & Compute $\fG_\textrm{X,s} := \rT\fP_s\rT\fX \fX \fP_s$
 \\
\hline
$[\fPhi, \fD] := \eig(\fG_\textrm{X,s})$ &
$[\fPhi, \fD] := \eig(\fG_\textrm{X,s})$\\
\hline
\multicolumn{2}{|c|}{\rule{0pt}{4ex}Set $ \fS := $diag$\left(\left(\frac{1}{\sqrt{\lambda_1}}, ... ,\frac{1}{\sqrt{\lambda_{m}}}\right)\right) \in \R^{n_s \times m},$}\\ \multicolumn{2}{|c|}{with $\lambda_1\geq... \geq\lambda_{m}>0$ the dominant, non-zero diagonal elements of $\fD$}  \\
\hline
Set $\widetilde \fPhi := \fPhi \fS$ &
Set $\fV :=\fX \fP_s \fPhi\fS $  \\
\hline
Select sub-matrix $\fAXs := \rT\fP_s\fAX\fP_s$  &  \\
\hline
Compute $ \fA_r := \rT{\widetilde \fPhi}\fAXs\widetilde \fPhi$&
Compute $\fA_r := \rT\fV\fA\fV$  \\
\hline
Select sub-vector $\xInitXs := \rT\fP_s\xInitX$ &\\
\hline
Compute $\fxrInit := \rT{\widetilde\fPhi}\xInitXs$& Compute $\fxrInit := \rT \fV\fxInit$\\
\hline
\end{tabular}}
\end{table}

The algorithms in both columns are identical in exact arithmetic, as the matrices are the same, but just the order of multiplications is changed. 
The matrix product $\fA_r =\rT{\widetilde \fPhi}\fAXs \widetilde \fPhi$ is online-efficiently computable as the operations depend only on the sizes of $
\fAXs \in \R^{\ns\times\ns}$ and $\widetilde\fPhi\in\R^{\ns \times m}, m \leq \ns$ %
i.e. the number of the online selected snapshots. 
After computation of a new basis, the reduced state needs to be projected to the new space, which also can be realized in an online-efficient manner:
Let $\widetilde\fPhi^i\in \R^{\ns\times m_i}$ denote the matrix of eigenvectors corresponding to the $i$th basis matrix, $\fP_s^i\in\R^{\NX\times\ns}$ the selection matrix for the $i$th set of indices and $\fxr^i$ a vector of reduced coefficients. 
If the basis is changed after the calculation of $\fxr^i$, the new representation $\fxr^{i+1} \in \R^{m_{i+1}}$ for the $(i+1)$th basis $\widetilde\fPhi^{i+1}\in \R^{\ns\times m_{i+1}}$ with selection matrix $\fP_s^{i+1}\in\R^{\NX\times\ns} $ can be calculated as
\begin{align*}\fxr^{i+1} = (\widetilde \fPhi^{i+1}\rT)(\fP_s^{i+1}\rT)\fGX\fP_s^i\widetilde \fPhi^i \fxr^i.
\end{align*}
\subsection{Dictionary-based DEIM}
\label{secDDEIM}
Next, we explain how an online-efficient, dictionary-based version of the DEIM-algorithm (DB-DEIM) can be realized. As a first step, we split the RHS of the FOM \Cref{eq:FOM}
$$\ff(\fx(t; \fmu), t;  \fmu) = \fA(t;\fmu)\fx(t; \fmu) + \ff_\textrm{nl}(\fx(t; \fmu),t; \fmu)$$
 into  a linear part $\fA(t;\fmu)$  and a non-linear part  $\ff_\textrm{nl}(\cdot,t; \fmu)$ and apply the DEIM-algorithm to the non-linear part. For notational simplicity we restrict to a non-parametric and time-independent $\fA, \ff_\textrm{nl}$. Parameter dependence can be straightforwardly added by considering all function evaluations to be parametric.
In order to efficiently treat non-linearities in the right-hand side of the ODE system in combination with the $\DEIMalg$-algorithm (see \Cref{DEIMalg}), a dictionary of non-linearity snapshots 
\begin{equation} 
  \dictF := \{\ff_{\text{nl}}(\xs_1), ..., \ff_{\text{nl}}(\xs_{\NX})\} \subset \R^{2N}
\end{equation}
and a dictionary of DEIM-indices $\dictP\subset \N$ is computed. The dictionary $\dictP$ is computed with $[\hat{\frho}^{\dictP}, \hat\fP_{\dictP}] := \DEIMalg(\hat \fU)$, where $\hat \fU$ = \texttt{POD}($\fF,\NP) \in \R^{2N \times \NP}$, $\NP \leq \NX$, with 
\[
  \fF:=[\ff_{\text{nl}}(\xs_1), ..., \ff_{\text{nl}}(\xs_{\NX})] \in \R^{2N \times \NX}.
\] 
The size of the dictionary of DEIM-indices $\NP \leq \NX$ has to be specified by the user. For simplicity, we use $\NP = \NX.$
Let $\hat\fP_{\dictP}=[\fe_{i_1}, ...,\fe_{i_{\NP}}] \in \R^{2N \times{\NP}}$, while $\fe_i$ being the $i$th unit vector, denote the selection matrix corresponding to the DEIM index-vector $\hat{\frho}^{\dictP} := (i_1, ..., i_{\NP}\rT)$. The dictionary $\dictP$ is chosen as the set of the entries of the DEIM index-vector $\hat{\frho}^{\dictP}$ 
\begin{equation}\label{eqndictP}
  \dictP := \{i_1, ..., i_{\NP}\} \subset \{1,...,2N\}.
\end{equation}

After the application of the standard DEIM-method, where the right-hand side $\ffr(\fxr(t; \fmu))$ is split into a linear and a non-linear part, the reduced ODE system (\ref{eq:ROM}) reads 
\begin{equation}
\ddt \fxr(t) = \rT\fV \fA \fV\fxr(t) + \rT\fV \fU (\rT\fP\fU) ^{-1} \fg(\fxr(t))
\end{equation}
with $\fg(\fy) = \rT\fP\ff_{\textrm{nl}}( \fV \fy), \fy\in \R^{m}$. 
In order to combine DEIM with dictionary-based MOR, an online-efficient computation of $\fGVU :=\rT\fV \fU,$ $\fUP :=\rT\fP \fU$ and $\rT\fP \ff_{\textrm{nl}}(\fV\fy)$  is required. In \Cref{tablestandvsdictbasedDEIM} we show how this is obtained with a dictionary-based approach based on the dictionaries introduced in the beginning of \Cref{sec:main} and what are the differences to the standard approach. 
The corresponding mathematical operations in both algorithms are again aligned with each other.
\begin{table}[h]
\caption{Comparison of dictionary-based DEIM vs.\ standard DEIM.}
\label{tablestandvsdictbasedDEIM}
\centering
{\renewcommand{\arraystretch}{\customstretch}%
\small
\begin{tabular}{|c|c|}
\hline
\bf Dictionary-based DEIM: & \bf Standard DEIM:\\
\hline
\bf Offline-DB-DEIM: &\\
\hline
Call Offline-DB-POD from \Cref{tablestandvsdictbasedPOD} & \\
Glue all $\fx(t_j, \fmu_i)\in \fD_{\text{X}}$ together as $\fX$ & \\
glue all $\ff_{\textrm{nl}}(\fx(t_j, \fmu_i))\in \fD_{\text{F}}$ together as $\fF$, &\\
Precompute $\fGXF := \rT\fX\fF, \fGF = \rT\fF\fF$ &\\
 $\hat{\fU} :=$ POD($\fF, \NP$), &\\ 
$[\hat{\frho}^{\dictP},  \hat \fP_{\dictP}] := $ \DEIMalg($\hat{\fU}$), $\fFhatP := \rT{\hat \fP_{\dictP}}\fF$ & \\
\hline
\bf Online-DB-DEIM: & \\
\hline
Compute $\widetilde \fPhi$ (Online-DB-POD) & Compute $\fV$:=POD($\fX\fP_s,m$)  \\
\hline 
Select sub-matrices $\fGXFs := \rT\fP_s\fGXF\fP_s,$&  Compute $\fGFs := \rT\fP_s\rT\fF\fF\fP_s$\\
$\fGFs := \rT\fP_s\fGF\fP_s, \fFhatPs :=  \fFhatP\fP_s$&  \\
\hline
$[\fPsi, \fD] := \eig(\fGFs)$& $[\fPsi, \fD] := \eig(\fGFs)$\\
\hline
\multicolumn{2}{|c|}{\rule{0pt}{4ex}Set $ \fS := $diag$\left(\left(\frac{1}{\sqrt{\lambda_1}}, ... ,\frac{1}{\sqrt{\lambda_{\widetilde m}}}\right)\right) \in \R^{n_s \times \widetilde m},$}\\ \multicolumn{2}{|c|}{with $\lambda_1\geq... \geq\lambda_{\widetilde m}>0$ the dominant, non-zero, diagonal elements of $\fD$}  \\
\hline
Set $\widetilde \fPsi := \fPsi \fS$&  Set $\fU :=\fF \fP_s \fPsi \fS$\\
\hline
Compute $\fGVU := \rT{\widetilde{\fPhi}}\fGXFs\widetilde \fPsi$&  Compute $\fGVU := \rT\fV\fU$\\
\hline
Compute $[\frho_{\textrm{o}}, \fP_{\textrm{o}}]:=$ \DEIMalg($\fFhatPs\widetilde \fPsi$) & Compute $[\frho, \fP]:= \DEIMalg(\fU)$\\
\hline
Select rows $\fUPPo := \rT\fP_{\textrm{o}}\fFhatPs\widetilde \fPsi$ & Select rows $\fUP:=\rT\fP\fU$\\
\hline
Comp. $\fg(\fy) := \rT\fP_{\textrm{o}}\rT{\hat\fP_{\dictP}}\ff_{\textrm{nl}}(\fX \fP_s \widetilde\fPhi \fy)$ & Comp. $\fg(\fy):=\rT\fP\ff_{\textrm{nl}}(\fV\fy)$ \\
\hline
\end{tabular}}
\end{table}
Again, the fact is used that the POD bases for $\fX\fP_s$ and $\fF\fP_s$ can be expressed by a product of $\fX\fP_s$ or $\fF\fP_s$ and the eigenvectors of the small eigenvalue problems $\rT\fP_s\rT\fX \fX \fP_s\in\R^{\ns\times\ns}$ or $\rT\fP_s\rT\fF\fF\fP_s\in\R^{\ns\times\ns}$ scaling by the inverses of the square roots of the corresponding eigenvalues. Therefore, the necessary multiplications with $\fX$ and $\fF$ are precomputed, the desired entries of the matrix products are selected and multiplied with the corresponding scaled eigenvectors as a low-dimensional operation just depending on the number of selected snapshots. Like this, the matrix product $ \fGVU=\rT\fV \fU$ is calculated. In order to realize an online-efficient computation for $\fUP = \rT\fP \fU$ a pre-selection of entries has to be performed because the computational costs of the \DEIMalg-algorithm depend on the size of the input matrix. Thus, we compute a pre-selection $\fFhatP := \rT{\hat \fP_{\dictP}}\fF$ with $[\hat{\frho}^{\dictP}, \hat\fP_{\dictP}] = \DEIMalg(\hat \fU$) during the offline-phase. During the online-phase we then apply the \DEIMalg-algorithm to the small matrix $\fFhatPs\widetilde \fPsi \in \R^{\NP \times \ns}$ to compute an online-DEIM selection matrix $\fP_{\textrm{o}}$ and an online-index vector $\frho_{\textrm{o}}$ with 
$$[\frho_{\textrm{o}}, \fP_{\textrm{o}}]:= \DEIMalg(\fFhatPs\widetilde \fPsi),$$
which will then be used for the selection. The selection is then formally performed with the selection matrix $\hat \fP_{\dictP}\fP_{\textrm{o}}$. The selection $(\hat \fP_{\dictP}\fP_{\textrm{o}}\rT)\fU$ is then online-efficiently performed using the pre-computed selection $\fFhatPs$ and the scaled eigenvectors collected in $\widetilde\fPsi$ as 
$$\fUPPo := \rT\fP_{\textrm{o}}\fFhatPs\widetilde \fPsi = (\hat \fP_{\dictP}\fP_{\textrm{o}}\rT)\fU.$$
The selection of entries of $\ff_{\textrm{nl}}$ is then also performed using $\hat \fP_{\dictP}\fP_{\textrm{o}}$ as selection matrix as $\rT\fP_{\textrm{o}}\rT{\hat\fP_{\dictP}}\ff_{\textrm{nl}}(\fX \fP_s \widetilde\fPhi \fy)$. Note that the reconstructions $\fX \fP_s \widetilde\fPhi \fy$ and $\fV\fy$, respectively, are actually not fully computed in this step. For DEIM we assume, that each component of $\ff_\textrm{nl}$ depends only on a few entries of $\fX \fP_s \widetilde\fPhi \fy$ and $\fV\fy$, respectively, (which is sometimes referred as $2N$-independent local DOF dependence, \cite{Drohmann2012}). For a discretized PDE, this assumption on $\ff_\textrm{nl}$ usually holds. These entries are then selected and computed without forming the high-dimensional product $\fV\fy$ or $\fX \fP_s \widetilde\fPhi \fy$. This selection process depends on the structure of $\ff_\textrm{nl}$. Thus, we demonstrate the selection process for a special function $\ff_\textrm{nl}$ in our numerical experiments in \Cref{sec:numerics}. 
The computation of the matrix $\fGVU$ in the left and right columns is identical in exact arithmetic, as only the order of multiplications in the algorithms on the left and right differ. The selection matrices $\rT{\fP}$ and $\rT{\fP_{\textrm{o}}}\rT{\hat\fP_{\dictP}}$ in the two algorithms differ in general, because in the dictionary-based DEIM algorithm, an offline-pre-selection is performed during the DEIM-index-dictionary computation and then a second DEIM index selection is performed during the online-simulation. With the standard approach, the DEIM indices would be calculated directly in one step. However,  with the following proposition a condition is proven that ensures equivalence of the two algorithms. 

\begin{proposition}[Equivalence of DB-DEIM index selection]\label{PropDEIM} 
\quad   \\
\noindent
Let $[\hat \fP_{\dictP}, $$\hat{\frho}^{\dictP}]$$ :=  $$\DEIMalg(\hat{\fU})$ be the selection matrix and index vector corresponding to the dictionary $\dictP$. Let $[\fP,  \frho]= \DEIMalg(\fU)$,  and $[\fP_\textrm{o}, \frho^\textrm{o}] :=  $\DEIMalg$(\fFhatPs\widetilde \fPsi)$. If the set of selected indices $M_\frho :=\{\rho_1, ..., \rho_{\widetilde{m}}\vert (\rho_1, ..., \rho_{\widetilde{m}}\rT)= \frho \}$ is contained in the dictionary $ M_\frho\subset \dictP$ then $\hat \fP_{\dictP}\fP_\textrm{o} = \fP$ and $\frho = \hat{\frho}^{\dictP}(\frho^\textrm{o}):=\rTb{ \rho^{\dictP}_{\rho^\textrm{o}_1}, ..., \rho^{\dictP}_{\rho^\textrm{o}_{\widetilde m}}}$.
\end{proposition}
\vspace{-0.3cm}
\begin{proof} 
We start by showing, that $\rT{\hat{\fP}_{\dictP}}\fU = \fFhatPs\widetilde \fPsi:$ With the definitions from \Cref{tablestandvsdictbasedDEIM} it follows, that $\rT{\hat{\fP}_{\dictP}}\fU = \rT{\hat{\fP}_{\dictP}}\fF \fP_s \fPsi \fS = \rT{\hat{\fP}_{\dictP}}\fF \widetilde \fPsi = \fFhatPs\widetilde \fPsi.$ In the following, we will use the representation $\rT{\hat{\fP}_{\dictP}}\fU$ instead of $\fFhatPs\widetilde \fPsi$. 

Now, the rest of the proof follows by induction. We first show, that $\rho_1 = \rho^{\dictP}(\rho^\textrm{o}_1)$. From the $\DEIMalg-$Algorithm it follows, that 
$$\rho_1 = \textrm{argmax}(\vert\fU(:,1)\vert)$$ 
and 
$$\rho^\textrm{o}_1 = \textrm{argmax}(\vert\rT{\hat \fP_{\dictP}}\fU(:,1)\vert).$$
If $M_\frho :=\{\rho_1, ..., \rho_{\widetilde m}\vert (\rho_1, ..., \rho_{\widetilde m}\rT) = \frho \}\subset \dictP$ then $$\max\vert\fU(:,1)\vert = \max\vert\rT{\hat \fP_{\dictP}}\fU(:,1)\vert$$ because $\vert\rT{\hat \fP_{\dictP}}U(:,1)\vert$ is a selection of $\vert\fU(:,1)\vert$ which contains the maximum entry of $\vert\fU(:,1)\vert$ as $\rho_1\in\dictP$. Thus, $\rho_1 = \rho^{\dictP}(\rho^\textrm{o}_1)$ and the first column of $\hat \fP_{\dictP}\fP_\textrm{o}$ and $\fP$ are identical. 

For the induction step we assume that the first $\ell-1, \ell\geq 2$ columns of $\hat \fP_{\dictP}\fP_\textrm{o}$ and $\fP$ are identical and $\rho_i = \rho^{\dictP}(\rho^\textrm{o}_i), i =1,...,\ell-1.$
The next step of the $\DEIMalg-$Algorithm is to calculate
\begin{equation} \label{cDEIMP}
\fc= (\fP(:,1:\ell-1\rT)\fU(:,1:\ell-1))^{-1}(\fP(:,1:\ell-1\rT)\fU(:,\ell)).
\end{equation} 
During the calculation of $[\fP_\textrm{o}, \frho^\textrm{o}] =  \DEIMalg(\rT{\hat{\fP}_{\dictP}}\hat{\fU})$ this leads to 
\begin{equation}\label{cDEIMPo}
\fc^\textrm{o}=(\fP_\textrm{o}(:,1:\ell-1\rT)\rT{\hat{\fP}_{\dictP}}\fU(:,1:\ell-1))^{-1}(\fP_\textrm{o}(:,1:\ell-1\rT)\rT{\hat{\fP}_{\dictP}}\fU(:,\ell)).
\end{equation}
As ${\hat{\fP}_{\dictP}}\fP_\textrm{o}(:,1:\ell-1)$ are the first $\ell-1$ columns of ${\hat{\fP}_{\dictP}}\fP_\textrm{o}$, and they are assumed to be identical to the first $\ell-1$ columns of $\fP$ by the induction assumption, the calculation in \Cref{cDEIMP,cDEIMPo} are identical and $\fc^\textrm{o} = \fc.$
From the $\DEIMalg-$Algorithm it then follows, that 
$$\rho_\ell = \textrm{argmax}(\vert\fU(:,\ell) - \fU(:,1:\ell-1)\fc\vert)$$
and
$$\rho^\textrm{o}_\ell = \textrm{argmax}(\vert\rT{\hat \fP_{\dictP}}(\fU(: - \fU(:,1:\ell-1))\fc^\textrm{o}\vert).$$
Because we assumed $\rho_\ell\in \dictP$, and it holds that $\fc^\textrm{o} = \fc$, it follows that the maximum entry of $\vert(\vert\fU(:,\ell) - \fU(:,1:\ell-1))\fc\vert$ will also be the maximum entry of the selection $\vert\rT{\hat \fP_{\dictP}}(\fU(: - \fU(:,1:\ell-1))\fc^\textrm{o})\vert.$ Thus, also the $l$th column of $\hat \fP_{\dictP}\fP_\textrm{o}$ and $\fP$ are identical and $\rho_\ell = \rho^{\dictP}(\rho^\textrm{o}_\ell).$ Thus, it follows by induction, that 
$$\hat \fP_{\dictP}\fP_\textrm{o} = \fP\ \ \textrm{and} \ \ \frho = \hat{\frho}^{\dictP}(\frho^\textrm{o})=(\rho^{\dictP}_{\rho^\textrm{o}_1}, ..., \rho^{\dictP}_{\rho^\textrm{o}_{\widetilde m}}\rT).$$
\end{proof}
We learn from Proposition \ref{PropDEIM} that if the dictionary contains a sufficient number of elements, the appropriate indices will be available and selected during the online-phase. 
\subsection{Dictionary-based Symplectic MOR}
\label{secDSMOR}
In the following, we discuss how the concepts of dictionary-based MOR
and symplectic MOR can be combined in order to formulate a structure-preserving basis computation (step~\circNum{2} and ~\circNum{3} of the dictionary-based basis update from \Cref{fig:Workflow}). We focus on a dictionary-based version of the cSVD-algorithm (DB-SVD). The cSVD basis from \cite{Peng2016} can either be computed by a singular value decomposition of a complex snapshot matrix or the singular value decomposition of a real, extended snapshot matrix $\fY := [\fXs, \JtN \fXs]$, which is proven to be equivalent in \cite{Buchfink2019}. 
For an SVD we assume the singular values to be sorted in decreasing order. That means that, as it is stated in \cite[Proposition~6]{Buchfink2019}, the columns have to be re-arranged in order to get a symplectic matrix. In Algorithm \ref{csvdpody} it is shown how this can be implemented by selecting every second vector of the POD, stacking them together for the first $k$ columns of $\fV $. 
\begin{algorithm}
\caption{Complex SVD via POD of $\fY$}
\label{csvdpody}
\begin{algorithmic}
\State Input: Snapshot matrix $\fXs\in \R^{2N\times n_{\fXs}}$, size $2n$ of the ROB matrix
\State Output: Symplectic basis $\fV \in \R^{2N\times 2n}$
\State $\fY = [\fXs, \JtN \fXs]$
\State $\fV_{\fY}$=  POD($\fY,2n$) \Comment{compute ROB matrix of size $2n$ with POD}
\State$\fV = [\fV_{\fY}(:, 1:2:2n), \TJtN \fV_{\fY}(:, 1:2:2n)]$ 
\end{algorithmic}
\end{algorithm}
For the columns $k+1$ to $2k$ of $\fV$ every second vector of the POD is extracted, and they are multiplied with $\TJtN$. This procedure works if the singular values of $\fY = [\fXs, \JtN \fXs]$ occur only in multiplicity of 2. In the very unlikely case that there are 4 or more singular values that are equal, a more careful selection is required. 
We call this version of the cSVD-algorithm \emph{cSVD via POD of $\fY$} and use it to derive a dictionary-based version, which is shown in \Cref{tablestandvsdictbasedcSVD}. For notional simplicity, we consider a linear, non-parametric Hamiltonian system, i.e. $\grad[\fx] \Ham(\fx) = \fH\fx$, $\fxInit(\fmu) = \fxInit$. Additional separable parameter dependence of $\fH$ and $\fxInit$ is feasable by performing the operations to the components of $\fH$ and $\fxInit$. Then, an additional online linear combination step is performed. 

\begin{table}[H]
\caption{Comparison of dictionary-based cSVD vs.\ standard cSVD via POD of $\fY$ for a linear non-parametric system.}
\label{tablestandvsdictbasedcSVD}
\centering
{\renewcommand{\arraystretch}{\customstretch}%
\small
\begin{tabular}{|c|c|}
\hline
\bf Dictionary-based cSVD: & \bf Stand. cSVD (of $\fX\fP_s$):\\
\hline
\bf Offline-DB-cSVD: &\\
\hline
Glue all $\fx(t_j, \fmu_i) \in \fD_{\text{X}}$ together as $\fX$,& \\ 
compute $\fGX := \rT\fX\fX$, $\fGXJ := \rT\fX\JtN\fX$,& \\
$ \fHX := \rT\fX \fH  \fX$, $\fHXJr := \rT\fX\fH\JtN\fX$,& \\
 $\fHXJl := \rT\fX \JtN \fH \fX$, & \\
$\fHXJJ := \rT\fX \JtN\fH\JtN \fX$,  & \\
$\fx_\textrm{0,X} := \rT\fX \fxInit$ and  & \\
$\fx_{\textrm{0,X},\JtN } := \rT\fX \JtN\fxInit$& \\
\hline
\bf Online-DB-cSVD: & \\
\hline
Select sub-matrices $\fGXs := \rT\fP_s\fGX\fP_s$ & Compute $\fGYs:= \rT\fY_s\fY_s$ \\
and $\fGXJs := \rT\fP_s\fGXJ\fP_s$ & $\textrm{with } \fY_s := \fY[\fP_s; \fP_s]$\\
\hline
\rule{0pt}{4ex}
$ \rT\fPhi \fD \fPhi := \eig\left(\begin{pmatrix}
\fGXs \ \ & \fGXJs \\
 -\fGXJs \ \  &\fGXs
\end{pmatrix} \right)$ & 
$\begin{matrix}
 \rT\fPhi \fD \fPhi := \eig\left(\fGYs\right)  \\
\end{matrix}$\\
\hline
\multicolumn{2}{|c|}{\rule{0pt}{4ex}Set $\widetilde \fPhi = \left[\frac{\fPhi_1}{\sqrt{\lambda_1}}, \frac{\fPhi_3}{\sqrt{\lambda_3}}, ... \frac{\fPhi_{2k-1}}{\sqrt{\lambda_{2k-1}}}\right] \in \R^{2n_s \times k},$ with $\fPhi = [\fPhi_1,...,\fPhi_{\ns}]$} \\ 
\multicolumn{2}{|c|}{and $\lambda_1\geq\lambda_2\geq ... \geq\lambda_{2k}>0$ the dominant, non-zero, diagonal elements of $\fD$}  \\
\hline Select sub-matrices $ \fHXs := \rT\fP_s\fHX\fP_s$,& \\
 $\fHXJrs := \rT\fP_s\fHXJr\fP_s$,& \\
 $\fHXJls := \rT\fP_s\fHXJl\fP_s$ & \\
and $\fHXJJs := \rT\fP_s\fHXJJ\fP_s$ & \\
\hline
\rule{0pt}{5ex}Compute $\fA_r  := \Jtn
\begin{pmatrix}
\rT{\widetilde \fPhi} \fM_1\widetilde \fPhi \ \ \ & \rT{\widetilde \fPhi}\fM_2\widetilde \fPhi \\
\rT {\widetilde \fPhi} \fM_3\widetilde \fPhi \ \ \ & \rT{\widetilde \fPhi}\fM_4\widetilde \fPhi 
\end{pmatrix} $  & 
$\begin{matrix} \textrm{Set } \fV :=[\fY_s\widetilde \fPhi ,\TJtN \fY_s\widetilde \fPhi] \\
 \textrm{with } \fY_s := \fY[\fP_s; \fP_s] 
\end{matrix}$ \\
 & Compute $\fA_r  := \Jtn\rT\fV\fH\fV$   \\
\hline
Select sub-vectors
$\fx_\textrm{0,X,\textrm{s}} := \rT\fP_s\fx_\textrm{0,X}$ and  & \\
$\fx_{\textrm{0,X},\JtN, \textrm{s}} := \rT\fP_s\fx_{\textrm{0,X},\JtN }$ \\
\hline
\rule{0pt}{5ex} Compute $\fxrInit := \Jtn
\begin{pmatrix}
\rT{\widetilde \fPhi}
\begin{pmatrix}
\fx_{\textrm{0,X},\JtN, \textrm{s}} \\
\fx_\textrm{0,X,\textrm{s}} \
\end{pmatrix}\\
\rT{\widetilde \fPhi}
\begin{pmatrix}
-\fx_\textrm{0,X,\textrm{s}}\\
\fx_{\textrm{0,X},\JtN, \textrm{s}}\\
\end{pmatrix}
\end{pmatrix}$&
Compute $\fxrInit := \si{\fV} \fxInit$ \\
\hline
\end{tabular}}
\end{table}  
The matrices $\fM_1$ to $\fM_4$ in \Cref{tablestandvsdictbasedcSVD} are defined as 
 \[ \fM_1 := 
\begin{pmatrix}
\fHXs \  & \fHXJrs \\
 -\fHXJls \  &-\fHXJJs
\end{pmatrix},
 \fM_2 := 
\begin{pmatrix}
-\fHXJrs \ &\fHXs \\
 \fHXJJs \  &-\fHXJls
\end{pmatrix}, \]
\[
\fM_3 := 
\begin{pmatrix}
\fHXJls \  &\fHXJJs \\
 \fHXs \  &\fHXJrs
\end{pmatrix} \textrm{and }
\fM_4 := 
\begin{pmatrix}
-\fHXJJs \ &\fHXJls \\
 -\fHXJrs \  &\fHXs
\end{pmatrix}.\]

Let $\fPhi^i\in \R^{2\ns\times 2k_i}$ denote the matrix of eigenvectors corresponding to the $i$th basis matrix, $\fP_s^i\in\R^{\NX\times\ns}$ the selection matrix for the $i$th set of indices and $\fxr^{i}$ a vector of reduced coefficients. 
If the basis is changed after the calculation of $\fxr^i \in \R^{2k}$  the new representation $\fxr^{i+1} \in \R^{2k_{i+1}}$ of $\fxr^i \in \R^{2k}$ can be calculated as symplectic projection
\begin{align*}
\fxr^{i+1} = \rT\Jtn (\fV^{i+1}\rT) \JtN \fV^i \fxr^i.
\end{align*}
For this projection after a basis change, the basis matrices $\fV^{i}$ and $\fV^{i+1}$ are expressed as
\[\fV^i = [\fY[\fP_s^i,\fP_s^i]\widetilde \fPhi^i , \rT \JtN \fY[\fP_s^i,\fP_s^i]\widetilde \fPhi^i] \] and
\[\fV^{i+1} = [\fY[\fP_s^{i+1},\fP_s^{i+1}]\widetilde \fPhi^{i+1} , \rT \JtN \fY[\fP_s^{i+1},\fP_s^{i+1}]\widetilde \fPhi^{i+1}] \]
with $ \fY = [\fX, \TJtN \fX]$ and $\widetilde \fPhi$ from \Cref{tablestandvsdictbasedcSVD}.
This is online-efficiently computable as
\begin{equation}
\label{projectionintonewspace}
(\fV^{i+1}\rT) \JtN (\fV^i)
 = 
\begin{pmatrix}
(\widetilde \fPhi^{i+1}\rT) \fM_5\widetilde \fPhi^i \ \ \ & (\widetilde \fPhi^{i+1}\rT)\fM_6\widetilde \fPhi^i \\
 -(\widetilde \fPhi^{i+1}\rT) \fM_6\widetilde \fPhi^i \ \ \ & (\widetilde \fPhi^{i+1}\rT) \fM_5\widetilde \fPhi^i
\end{pmatrix}
\end{equation}

with 
\begin{align*}
\fM_5 :=
\begin{pmatrix}
(\fP_s^{i+1}\rT)\fGXJ\fP_s^i \ \ \ &-(\fP_s^{i+1}\rT)\fGX\fP_s^i \\
(\fP_s^{i+1}\rT) \fGX\fP_s^i \  &(\fP_s^{i+1}\rT)\fGXJ\fP_s^i
\end{pmatrix}, 
\end{align*}

\begin{align*}
\fM_6 :=
\begin{pmatrix}
(\fP_s^{i+1}\rT)\fGX\fP_s^i \ \ \ &(\fP_s^{i+1}\rT)\fGXJ\fP_s^i \\
 -(\fP_s^{i+1}\rT)\fGXJ\fP_s^i \  &(\fP_s^{i+1}\rT)\fGX\fP_s^i
\end{pmatrix}. 
\end{align*}
For the SDEIM-algorithm \cite{Peng2016}, a dictionary-based online-efficient version can be derived similarly to the dictionary-based DEIM-algorithm. 
We split the gradient of the Hamiltonian 
$$\grad[\fx] \Ham(\fx(t; \fmu); \fmu) = \fH(\fmu)\fx(t; \fmu) + \ff_\textrm{nl}(\fx(t; \fmu);\fmu)$$
into  a linear part $\fH(\fmu)$  and a non-linear part  $\ff_\textrm{nl}(\cdot; \fmu)$, which corresponds to a split in the RHS of the FOM \Cref{eq:FOM} 
$$\ff(\fx(t; \fmu),t;\fmu) =\JtN \fH(\fmu) \fx(t; \fmu) + \JtN \ff_\textrm{nl}(\fx(t; \fmu);\fmu)$$ 
and apply the DEIM-algorithm to the non-linear part. For notational simplicity we again restrict to a non-parametric $\fH, \ff_\textrm{nl}$. Similar to the DB-DEIM for DB-SDEIM we need a dictionary of non-linearity snapshots 
\begin{equation} 
  \dictF := \{\ff_{\text{nl}}(\xs_1), ..., \ff_{\text{nl}}(\xs_{\NX})\} \subset \R^{2N}
\end{equation}
and a dictionary of DEIM-indices $\dictP\subset \N$. With the definitions
\[
  \fF:=[\ff_{\text{nl}}(\xs_1), ..., \ff_{\text{nl}}(\xs_{\NX})] \in \R^{2N \times \NX}
\] and
$\hat \fU$ := \texttt{POD}($\fF,\NP) \in \R^{2N \times \NP}$, $\NP \leq \NX$
the dictionary $\dictP$ is computed with $[\hat{\frho}^{\dictP}, \hat\fP_{\dictP}] := \DEIMalg(\hat \fU$). We chose the dictionary $\dictP$ as the set of the entries of the DEIM index-vector $\hat{\frho}^{\dictP}$ 
\begin{equation}
  \dictP := \{i_1, ..., i_{\NP}\} \subset \{1,...,2N\}.
\end{equation}
After the application of the SDEIM-method, the reduced ODE system (\ref{eq:ROM}) reads 
\begin{equation}\label{eqn:SDEIM}
  \ddt { \fxr(t)} = \Jtn \rT\fV\fH\fV \fxr(t) + \Jtn \rT \fV\fU (\rT \fP\fU) ^{-1} \fg(\fxr(t))
\end{equation}
with $\fg(\fy) = \rT\fP\ff_{\textrm{nl}}(\fV \fy), \fy\in \R^{m}$.
Note, that because of the insertion of the term $\fU (\rT\fP\fU) ^{-1} \rT\fP$ the system in general can no longer be written as a Hamiltonian system since $\rT\fV\fU (\rT \fP\fU) ^{-1} \rT\fP\ff_{\textrm{nl}}(\fV \fxr)$ in general can not be expressed as the gradient of a potential $\Ham_\textrm{SDEIM}$. But if the SDEIM yields a good approximation, solving the equation should not lead to a large energy variation (see \cite{Peng2016}).
Because of the structure of \Cref{eqn:SDEIM}, an online-efficient computation of $\fGVU:= \rT\fV \fU, \fUP:=\rT\fP \fU$ and $\rT\fP\ff_{\textrm{nl}}(\fV \fy)$ for an arbitrary vector $\fy\in \R^{2k}$ is required. In \Cref{tablestandvsdictbasedSDEIM} is shown how this is obtained and what the differences to the standard approach are. The corresponding mathematical operations in both algorithms are again aligned with each other. The procedure is similar to the DB-DEIM, to compute a DB-SDEIM approximation also the selection $\fUP:=\rT\fP \fU$  has to be split in an offline-selection $\fFhatP := \rT{\hat \fP_{\dictP}}\fF$ and an online-selection. The formal online-selection $(\hat \fP_{\dictP}\fP_{\textrm{o}}\rT)\fU$ is online-efficiently computed using $\fFhatPs$ and $\widetilde \fPsi$ as
 $$\fUPPo :=\rT\fP_{\textrm{o}}\fFhatPs\widetilde \fPsi=(\hat \fP_{\dictP}\fP_{\textrm{o}}\rT)\fU.$$

\begin{table}[H]
\caption{Comparison of dictionary-based SDEIM vs.\ standard SDEIM}
\label{tablestandvsdictbasedSDEIM}
\centering
{\renewcommand{\arraystretch}{\customstretch}%
\small
\begin{tabular}{|c|c|}
\hline
\bf Dictionary-based SDEIM: & \bf Standard SDEIM:\\
\hline
\bf Offline-DB-SDEIM: &\\
\hline
Call Offline-DB-cSVD&\\
Glue all $\fx(t_j, \fmu_i)\in \fD_{\text{X}}$ together as $\fX$ & \\
glue all $\ff_{\textrm{nl}}(\fx(t_j, \fmu_i))\in \fD_{\text{F}}$ together as $\fF$, &\\
Precompute $\fGXF := \rT\fX\fF, \fGF := \rT\fF\fF$,  &\\
$\fGXFJ := \rT\fX\JtN\fF,\hat{\fU} :=$ POD($\fF, \NP$), &\\ 
$[\hat{\frho}^{\dictP},  \hat \fP_{\dictP}] := $ \DEIMalg($\hat{\fU}$), $\fFhatP := \rT{\hat \fP_{\dictP}}\fF$ & \\
\hline
\bf Online-DB-SDEIM: & \\
\hline
Compute $\widetilde \fPhi$ (Online-DB-cSVD) & Compute $\fV$:=cSVD($\fX\fP_s$)  \\
\hline 
Select sub-matrices $\fGXFs := \rT\fP_s\fGXF\fP_s,$& Compute $\fGFs :=\rT\fP_s\rT\fF\fF\fP_s$ \\
$\fGFs := \rT\fP_s\fGF\fP_s, \fFhatPs :=  \fFhatP\fP_s$&  \\
\hline
$[\fPsi \fD]:= \eig(\fGFs)$& $[\fPsi \fD]:= \eig(\fGFs)$\\
\hline
\multicolumn{2}{|c|}{\rule{0pt}{4ex}Set $ \fS := $diag$\left(\left(\frac{1}{\sqrt{\lambda_1}}, ... ,\frac{1}{\sqrt{\lambda_{\widetilde m}}}\right)\right) \in \R^{n_s \times \widetilde m},$}\\ \multicolumn{2}{|c|}{with $\lambda_1\geq... \geq\lambda_{\widetilde m}>0$ the dominant, non-zero, diagonal elements of $\fD$}  \\
\hline
Set $\widetilde \fPsi := \fPsi \fS$&  Set $\fU :=\fF\fP_s \fPsi \fS$\\
\hline
Select $\fGXFs := \rT\fP_s\fGXF\fP_s$ &\\
and  $\fGXFJs := \rT\fP_s\fGXFJ\fP_s$, &\\
\hline
\rule{0pt}{5ex}
Compute  $ \fGVU := \begin{pmatrix}
\rT{\widetilde \fPhi} \fM_7\widetilde \fPsi \\ 
\rT{\widetilde \fPhi} \fM_8\widetilde \fPsi 
\end{pmatrix} $&  Compute $\fGVU := \rT\fV\fU$\\
\hline
Compute $[\frho_{\textrm{o}},\fP_{\textrm{o}}]:=$ \DEIMalg($\fFhatPs\widetilde \fPsi$) & Compute $\fP:=$ \DEIMalg($\fU$)\\
\hline
Select rows $\fUPPo := \rT\fP_{\textrm{o}}\fFhatPs\widetilde \fPsi$ & Select rows $\fUP:=\rT\fP\fU$\\
\hline
Comp. $\fg(\fy) := \rT\fP_{\textrm{o}}\rT{\hat\fP_{\dictP}}\ff_{\textrm{nl}}([\fY_s\widetilde \fPhi ,\TJtN \fY_s\widetilde \fPhi] \fy)$ & Comp. $\fg(\fy):=\rT\fP\ff_{\textrm{nl}}(\fV\fy)$\\
with $\fY_s := [\fX, \JtN\fX][\fP_s; \fP_s] $  &\\
\hline
\end{tabular}}
\end{table} 
The matrices $\fM_7$ and $\fM_8$ from Table \ref{tablestandvsdictbasedSDEIM} are defined as
\[\fM_7 := 
\begin{pmatrix}
\fGXFs  \\
\fGXFJs \  
\end{pmatrix} \textrm{and }
\fM_8 := 
\begin{pmatrix}
\fGXFJs \\
\fGXFs 
\end{pmatrix}.\]
As for the DB-DEIM, the basis matrix $\fV$ in these steps is not actually completely computed,  but only the required entries of $\fV\fy$ are extracted. 

Now, that all method ingredients have been introduced, we want to summarize the procedure: The algorithm starts with the selection of snapshots from the dictionary with \Cref{Selalg} and start time $\tInit$. Next, the cSVD- and SDEIM approximations are computed using the DB-cSVD from \Cref{tablestandvsdictbasedcSVD} and DB-SDEIM from \Cref{tablestandvsdictbasedSDEIM}. Then, the start vector is projected using DB-cSVD. After that, the time stepping for the reduced system is advanced until a basis update is triggered. After a basis update has been executed, the procedure is repeated until we reach our final simulation time $\tEnd$. Instead of the start vector projection, Equation \eqref{projectionintonewspace} is used to find a representation of the current iterate in the next subspace. This dictionary-based (symplectic) MOR-algorithm can be similarly applied using the non-symplectic dictionary-based methods and equations from \Cref{secDPOD,secDDEIM}.

\subsection{Error Analysis}
In this section, we analyze, how the error in the Hamiltonian behaves during the application of our dictionary-based algorithm. Because of the symplectic basis generation the reduced Hamiltonian will be piecewise constant, if a symplectic integrator is used, that preserves the Hamiltonian, since the symplectic ROM conserves the Hamiltonian over time between the basis updates. During the projection from the old into the new subspace, when a basis update is triggered, an error in the Hamiltonian is introduced in general. This error can be bounded by the following proposition.   
\begin{proposition}[Basis Change Hamiltonian Error Bound]\label{lemmaprojerrorHam}\label{Hamerr}
Let $\fxr^i$ be the current iterate and $\fV^{i}$ and $\fV^{i+1}$ be the $i$th and $(i+1)$th basis matrices. Let $\fxr^{i+1} = (\fV^{i+1})^+ \fV^{i}\fxr^{i}$ and $\Ham \in C^1(\R^{2N}, \R)$ with $\|\grad[\fx] \Ham\|_\infty:= \sup\limits_{\fx \in \R^{2N}}\|\grad[\fx] \Ham(\fx)\|_2$ bounded. Then, 
\[ \vert\Ham(\fV^{i}\fxr^{i}) -  \Ham(\fV^{i+1}\fxr^{i+1})\vert \leq \|\grad[\fx] \Ham\|_\infty\|(\I{2N} - \fV^{i+1}(\fV^{i+1})^+)\fV^{i}\fxr^{i}\|_2.\]
\end{proposition}
\vspace{-0.3cm}
\begin{proof} 
By the mean value theorem there exists $\fx$ from the convex hull $\textrm{conv}(\fV^{i}\fxr^{i}, \fV^{i+1}\fxr^{i+1})$ such that 
\begin{align*}
\vert\Ham(\fV^{i}\fxr^{i}) -  \Ham(\fV^{i+1}\fxr^{i+1})\vert &= \vert\grad[\fx] \Ham(\fx\rT)(\fV^{i}\fxr^{i} - \fV^{i+1}\fxr^{i+1})\vert \\
 &\leq \|\grad[\fx] \Ham(\fx)\|_2\ \|\fV^{i}\fxr^{i} - \fV^{i+1}\fxr^{i+1}\|_2 \\
 &\leq \|\grad[\fx] \Ham\|_\infty\|\fV^{i}\fxr^{i} - \fV^{i+1}\fxr^{i+1}\|_2\\
& = \|\grad[\fx] \Ham\|_\infty\|(\I{2N} - \fV^{i+1}(\fV^{i+1})^+)\fV^{i}\fxr^{i}\|_2.
\end{align*}
\end{proof}
\begin{remark}
From \Cref{Hamerr} follows with the triangle inequality, that the error in the Hamiltonian $\vert\Ham(\fxInit(\fmu)) -\Ham(\fV\fxr(t, \mu)\vert$ is bounded by $\|\grad[\fx] \Ham\|_\infty$ times the sum of the projection errors plus a possible error from time discretization, if the integrator does not exactly conserve the Hamiltonian.

\end{remark}

\section{Numerical Experiments}\label{sec:numerics}
In this section, the developed online-adaptive, structure-preserving methods are applied to a linear and to a non-linear wave-equation model. We first present the 2D linear wave equation model. 
\subsection{2D linear wave equation}
\label{sec:2Dlinwav}
\begin{figure}[H]
\caption{Linear wave equation: Average relative reduction error over number of basis vectors and online-runtime, reproduction experiment}
\label{fig:methcomp1}
\begin{minipage}[t]{0.475\textwidth}
\includegraphics[width = \textwidth]{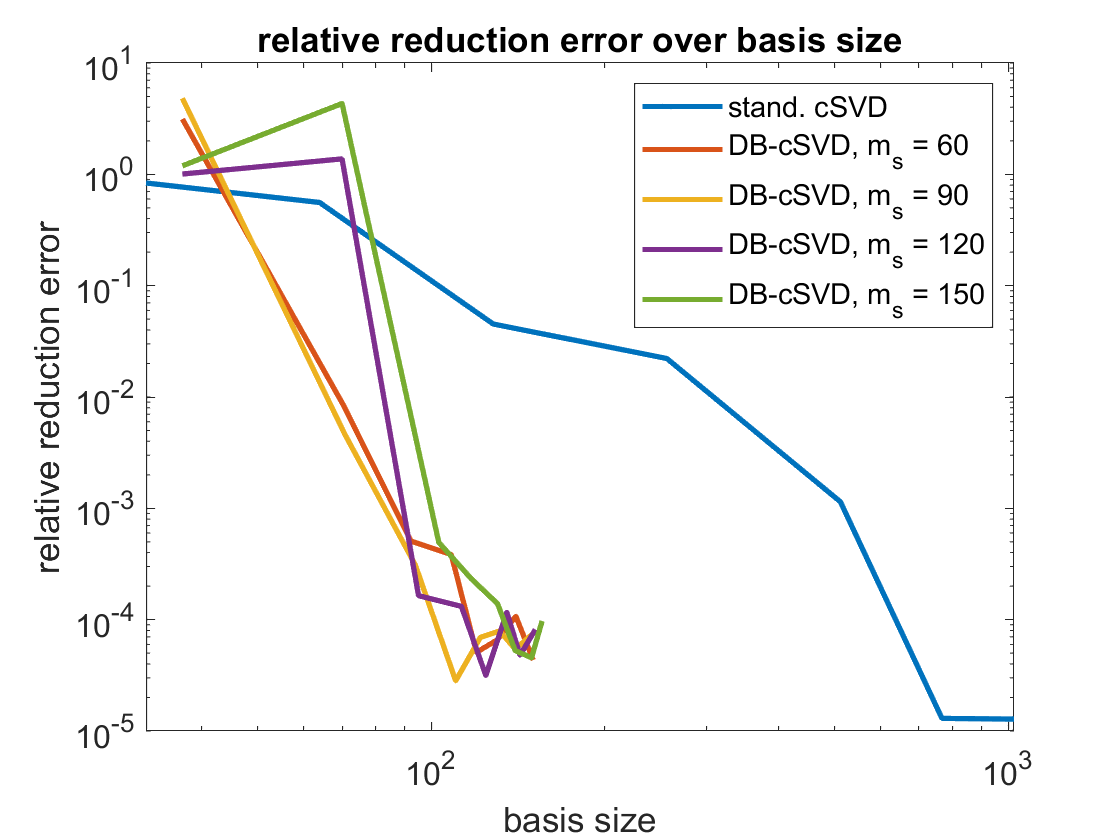}
\end{minipage}
\begin{minipage}[t]{0.475\textwidth}
\includegraphics[width = \textwidth]{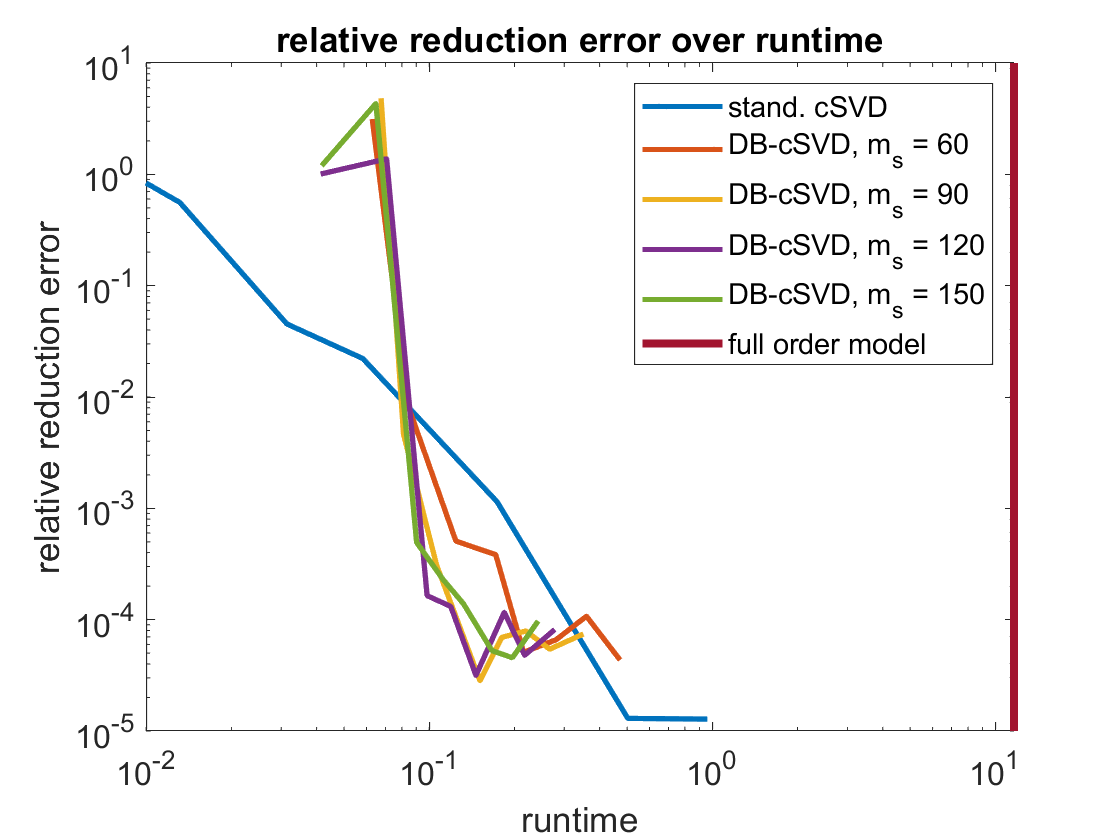}
\end{minipage}
\end{figure}

The initial boundary value problem for the unknown $u(t, \fxi)$ with the spatial variable $\fxi := (\xi_1, \xi_2) \in \varOmega := (0, l) \times (0, l/5), l:=1$ and the temporal variable $ t \in \It(\fmu) := [\tInit, \tEnd(\fmu)]$, reads
\begin{align*}
u_{tt}(t, \fxi) &= c^2 \Delta u(t, \fxi) &&\textrm{in } \It(\fmu)\times\varOmega\\
u(t_0, \fxi) &= u^0(\fxi) := h(s(\fxi)), \ u_t (t_0, \fxi) = v^0(\fxi) := 50c \ d_h(s(\fxi)) &&\textrm{in } \varOmega,\\
u(t, \fxi) &= 0 &&\textrm{in }  \It(\fmu)\times\partial\varOmega, 
\end{align*}
where
\[s(\fxi): = 50\cdot\left(\xi_1 - \frac{9l}{10}\right), \ 
h (s): = \begin{cases}
  1 - \frac{3}{2} \vert s\vert^2 + \frac{3}{4} \vert s\vert^3,
& 0 \leq \vert s\vert \leq 1, \\
  \frac{1}{4} (2 - \vert s\vert)^3,
& 1 < \vert s\vert \leq 2,  \\
  0,
& \vert s\vert > 2,
\end{cases}  \]
and
\[d_h (s): =
\begin{cases}
  - 3 s + \frac{9}{4} s^2,
& 0 \leq \vert s\vert \leq 1, \\
  \frac{1}{4} (-12\cdot\mathrm{sign(s)} + 12\vert s\vert - 3s\cdot\vert s\vert),
& 1 < \vert s\vert \leq 2,  \\
  0,
& \vert s\vert > 2.
\end{cases}  \]
We choose $\tInit = 0,\, \tEnd(\fmu) = 2/\fmu$ and as parameter (vector) $\fmu = c \in\mathcal{P} := [7,10]$.
The linear wave equation is spatially discretized using central finite differences.  This finally leads to the Hamiltonian system
\begin{equation}\label{wave_discr}
\ddt \fx(t; \fmu)  = \JtN \grad[\fx] \Ham(\fx(t; \fmu); \fmu) = \JtN\fH\fx, \quad \fx(0; \fmu) = \fxInit(\fmu),
\end{equation}
with 
$$\fH(\fmu) =\begin{pmatrix}
\fmu^2(\fD_{{\xi_1}{\xi_1}}+\fD_{{\xi_2}{\xi_2}}) \ &\Z{N} \\
 \Z{N}\  & \I{N}
\end{pmatrix}$$ 
and
$$\fxInit(\fmu) = [u^0(\fxi_1)), ...,u^0(\fxi_N)),v^0(\fxi_1)), ...,v^0(\fxi_N))]$$
where $\{ \fxi_i \}_{i=1}^N \subset \varOmega$ are the grid points.
The positive definite matrices $\fD_{{\xi_1}{\xi_1}}\in \R^{N \times N}$ and $\fD_{{\xi_2}{\xi_2}} \in \R^{N \times N}$ denote the three-point central difference approximations in $\xi_1-$direction and in $\xi_2-$direction. The number of grid points in $\xi_1-$ and $\xi_2$-direction are chosen as $N_{\xi_1} = 2000$ and $N_{\xi_2} = 20$, which results in a Hamiltonian system of dimension of $2N = 80000$. The corresponding Hamiltonian reads 
$$\Ham(\fx, \fmu) = \frac{1}{2}\rT\fx\fH(\fmu)\fx. $$ \\
Temporal discretization is achieved using the implicit midpoint rule  and $n_t = 600 $ equidistant time steps. Note that this leads to different time step sizes for different parameters because of the parameter dependence of $\tEnd(\fmu)$. The implicit midpoint belongs to the class of symplectic integrators \cite{hairer2006}, which is used for structure-preserving integration. Note however for later, that the implicit midpoint rule only preserves quadratic Hamiltonian functions exactly. 
First, we study a reproduction experiment. For the state dictionary we use the snapshots from the parameters, $\fmu_1 = 7, \fmu_2 = 8.5, \fmu_3 = 10$, from which also the snapshot matrix for the standard cSVD-basis is constructed. Thus, our dictionary size is $\NX = 1800.$
In \Cref{fig:methcomp1}, we present the relative reduction error 
\begin{equation}\label{eqnrelerr}
\fe_\text{rel}(\fmu) = \sqrt{\sum\limits_{i = 0}^{n_t}\|{\fx}_i(\fmu)  - \fV{\fxr}_i(\fmu)\|_2^2}\Bigg{/}\sqrt{\sum\limits_{i = 0}^{n_t}\|{\fx}_i(\fmu)\|_2^2}, 
\end{equation} 
with $  {\fx}_i(\fmu), {\fxr}_i(\fmu), i = 0,.., n_t$ the iterates of the FOM and ROM, in dependence on the average basis size and on the runtime of the reduced simulation (online-runtime). All reduction errors and runtimes are averaged over the three training parameters. For the DB-cSVD on the horizontal axis, the average number of the used basis vectors is computed as 
\begin{equation}\label{avgbasissisze}
n_\textrm{mean} = \sum\limits_{ j = 0}^{\lfloor\frac{n_t}{m_s}\rfloor+1} 2k_i.
\end{equation} 
The numbers $2k_i$ denote the basis size of the $i$th basis.  In each basis computation, the basis size is chosen as the lowest number $k \in \N$ that fulfills 
\begin{equation}\label{DBcSVDselcrit}
(1-\epsilon_\textrm{cSVD})\sum\limits_{\ell = 0}^{n_s}\lambda_{2\ell+1}^i<\sum\limits_{\ell = 0}^{k+1}\lambda_{2\ell+1}^i,
\end{equation}
where we denote with $\lambda_{\ell}^i>0, \ell = 1,..., 2k$ the dominant non-zero eigenvalues of 
$$\begin{pmatrix}
\fGXs^i \ \ & \fGXJs^i \\
 -\fGXJs^i \ \  &-\fGXs^i
\end{pmatrix},$$
with $\fGXs^i:=(\fP_s^i\rT) \fGX\fP_s^i$ and $\fGXJs^i :=(\fP_s^i\rT)\fGXJ \fP_s^i$, with $\fP_s^i$ the $i$th selection matrix and $\fGX$ from \Cref{tablestandvsdictbasedcSVD}. The eigenvalues occur in pairs, i.e. $\lambda_{2l+1}^i =\lambda_{2l+2}^i, l = 0,...,n_s.$ For the wave-equation experiments we chose $\epsilon_\textrm{cSVD}=10^{-12}.$ We compare the standard cSVD with the newly introduced DB-cSVD for different window sizes $m_s \in \{60, 90, 120, 150 \}$ and numbers of selected snapshots $\ns \in \{50, 100,...,400\}$ and the standard POD.
We observe that by using a dictionary-based approach, the same reduction error is achieved by far less basis vectors. For example a relative error of  $10^{-4}$ is achieved with an average basis size of about 100 with the dictionary-based methods compared to 500 basis vectors required for a standard cSVD basis. This reduced number of basis vectors also results in better runtimes. A relative error of $10^{-4}$ can be achieved with about 3 times shorter runtimes compared to the standard approach and nearly 2 orders of magnitude shorter compared to the FOM. The temporary increase in the error curves for $m_s \geq 120$ and an average basis size of 70 is due to the fact that the basis is too small for the window size $m_s\geq120$. Larger window sizes require larger basis sizes.

\begin{figure}[H]
\caption{Linear wave equation: Average relative error in Hamiltonian over time-steps, reproduction experiment}
\label{fig:methcomp1ham_non_adapt}
\begin{minipage}[t]{0.475\textwidth}
\includegraphics[width = \textwidth]{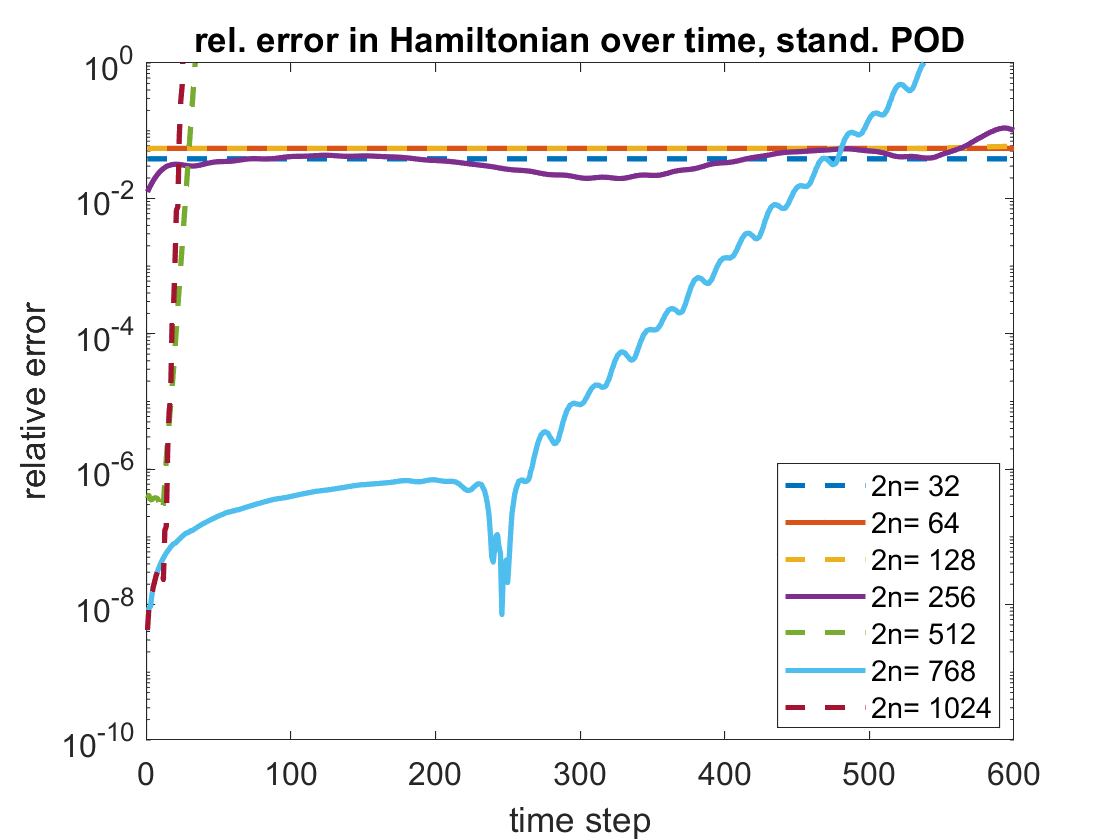}
\end{minipage}
\begin{minipage}[t]{0.475\textwidth}
\includegraphics[width = \textwidth]{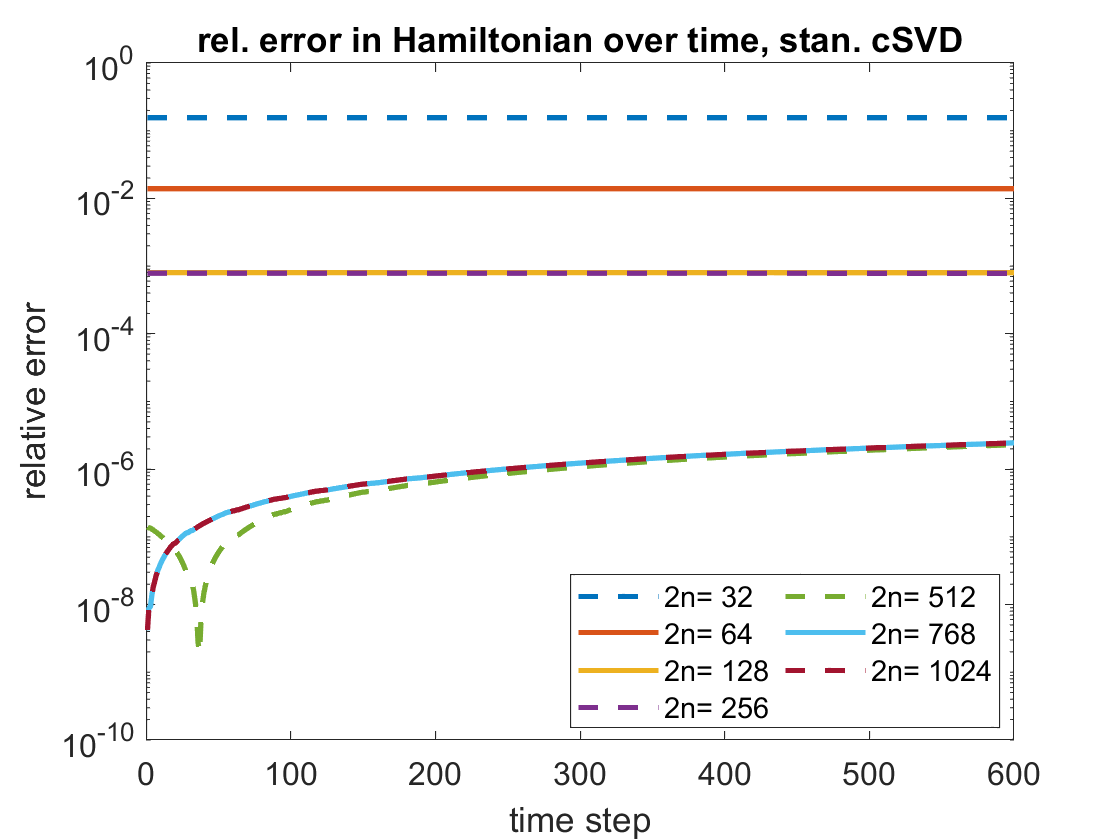}
\end{minipage}
\end{figure}

\begin{figure}[H]
\caption{Linear wave equation: Average relative error in Hamiltonian over time-steps, reproduction experiment}
\label{fig:methcomp1ham}
\begin{minipage}[t]{0.475\textwidth}
\includegraphics[width = \textwidth]{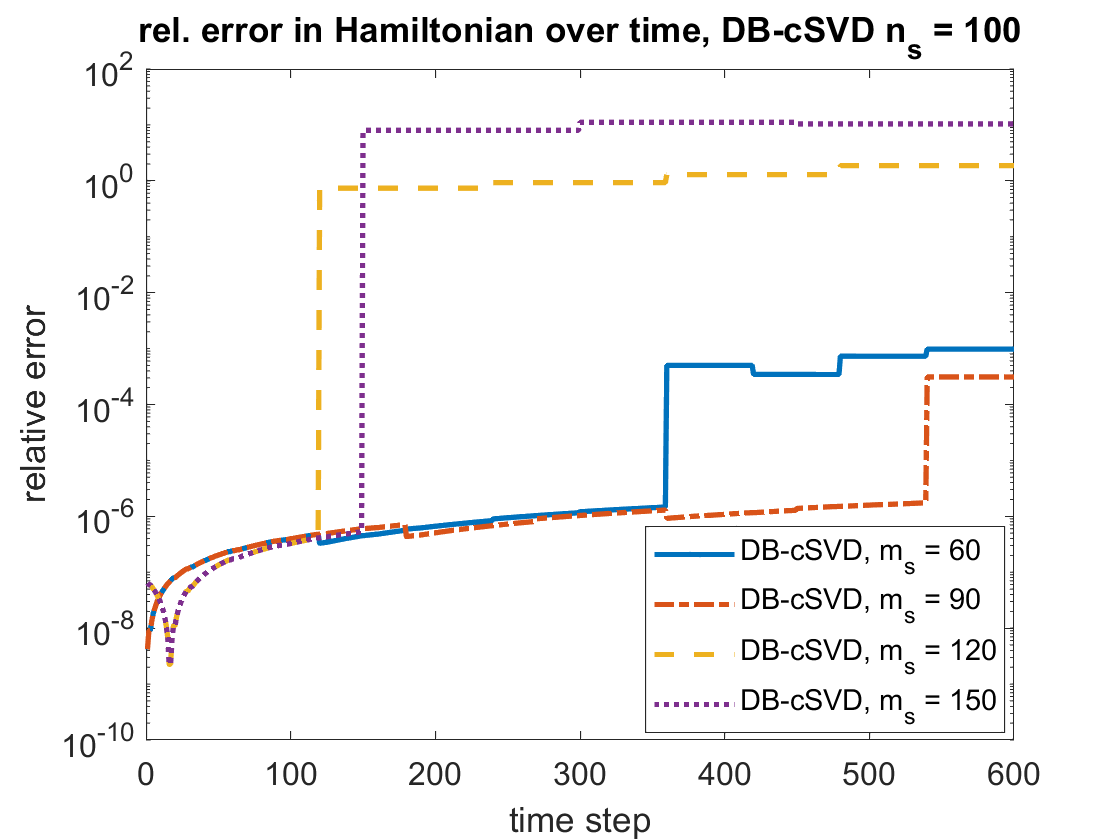}
\end{minipage}
\begin{minipage}[t]{0.475\textwidth}
\includegraphics[width = \textwidth]{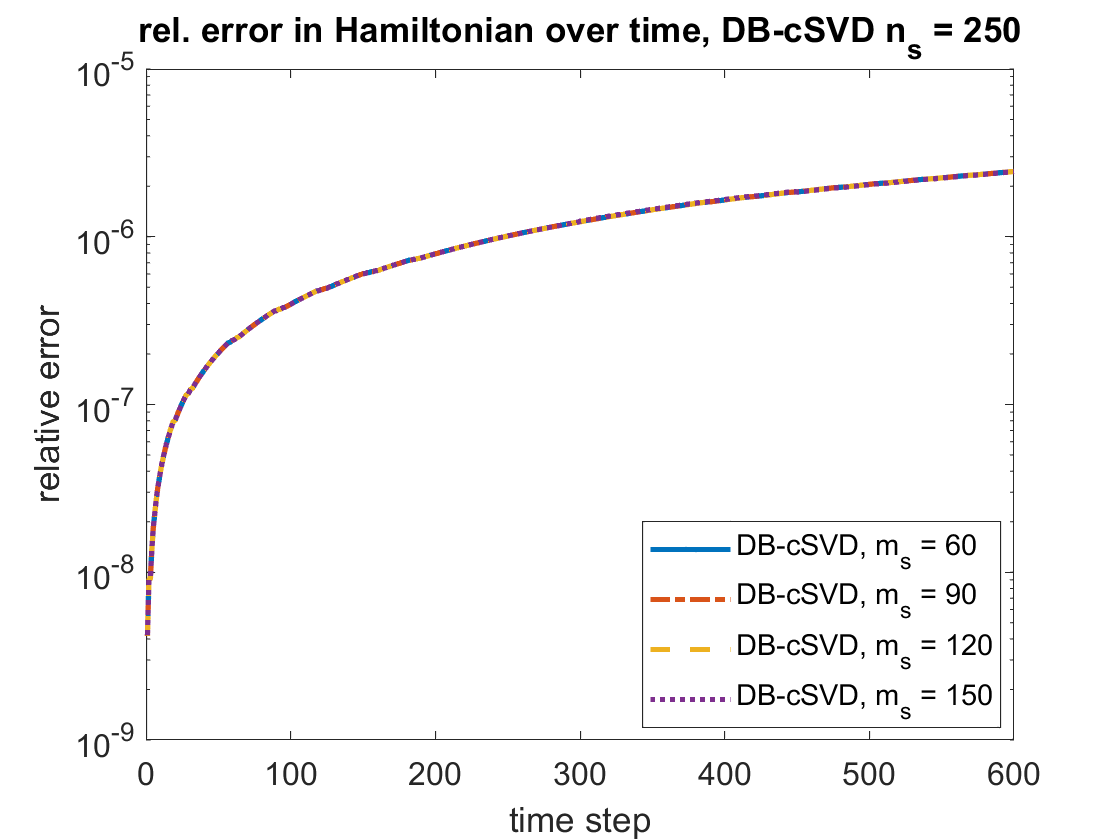}
\end{minipage}
\end{figure}

In \Cref{fig:methcomp1ham,fig:methcomp1ham_non_adapt} we present the relative error in the Hamiltonian 
\begin{equation}
\label{eqn:relHam}
\fe_{\Ham, \text{rel}, i}(\fmu) = \vert\Ham({\fx}_i(\fmu))  - \Ham(\fV{\fxr}_i(\fmu))\vert/\Ham({\fx}_i(\fmu)), 
\end{equation}
averaged over the three training parameters in dependence on the time-step.

In \Cref{fig:methcomp1ham_non_adapt} we see that using a non-symplectic, standard POD-basis leads to unstable models in the sense that the energy drastically increases if more than 256 basis vectors are used. 
In contrast, with 512 or more basis vectors the relative error in the Hamiltonian is about $10^{-6}$ for the standard cSVD which is in agreement with the experiments from \cite{Peng2016} for the linear wave equation. As expected, the error in the Hamiltonian is constant for the standard cSVD up to numerical inaccuracies. These inaccuracies can only be seen in the bottom two curves due to the logarithmic axis. In \Cref{fig:methcomp1ham}, we plot the same errors for the dictionary-based methods.
With $n_s = 100$ selected snapshots per basis update, we observe high error jumps across the basis changes for DB-cSVD for window sizes $m_s\geq 120$. However, with $n_s = 250$ snapshots selected, this behavior is no longer observed and a relative error of $10^{-6}$ is achieved as with the standard cSVD. 

\begin{figure}[H]
\caption{Linear wave equation: Average relative reduction error over number of basis vectors and online-runtime, generalization experiment}
\label{fig:methcomp1gen}
\begin{minipage}[t]{0.475\textwidth}
\includegraphics[width = \textwidth]{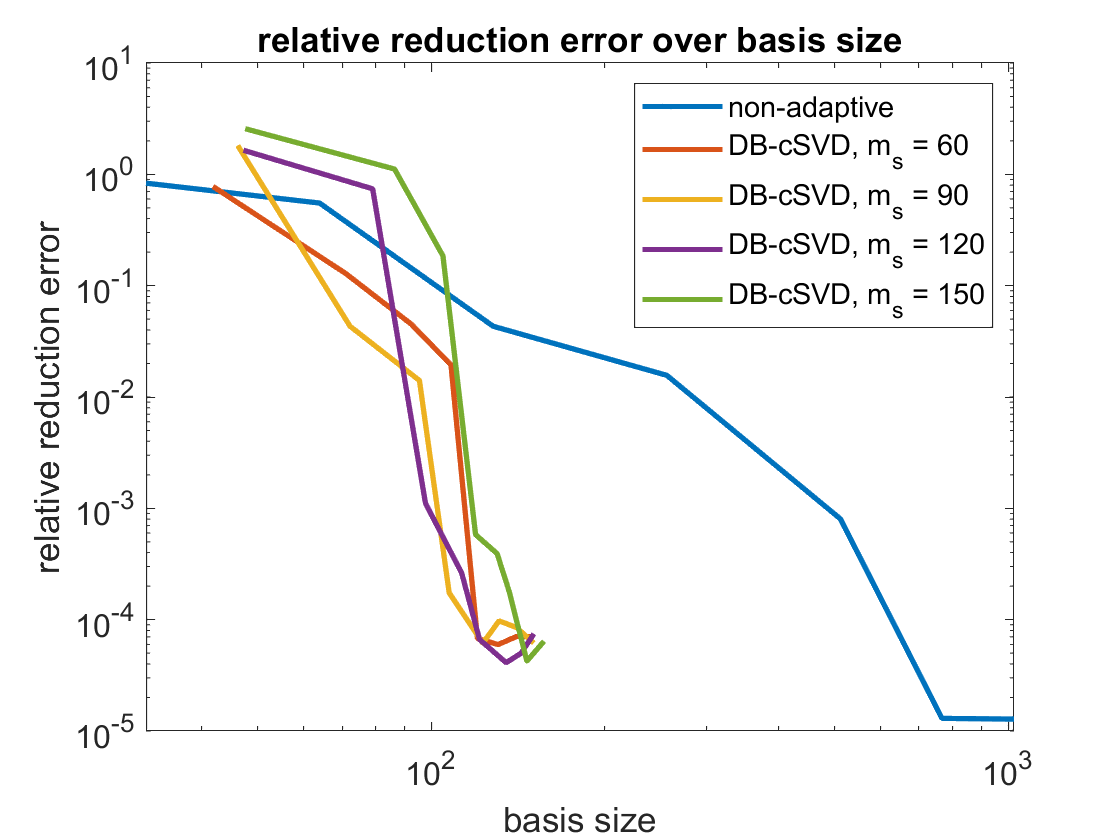}
\end{minipage}
\begin{minipage}[t]{0.475\textwidth}
\includegraphics[width = \textwidth]{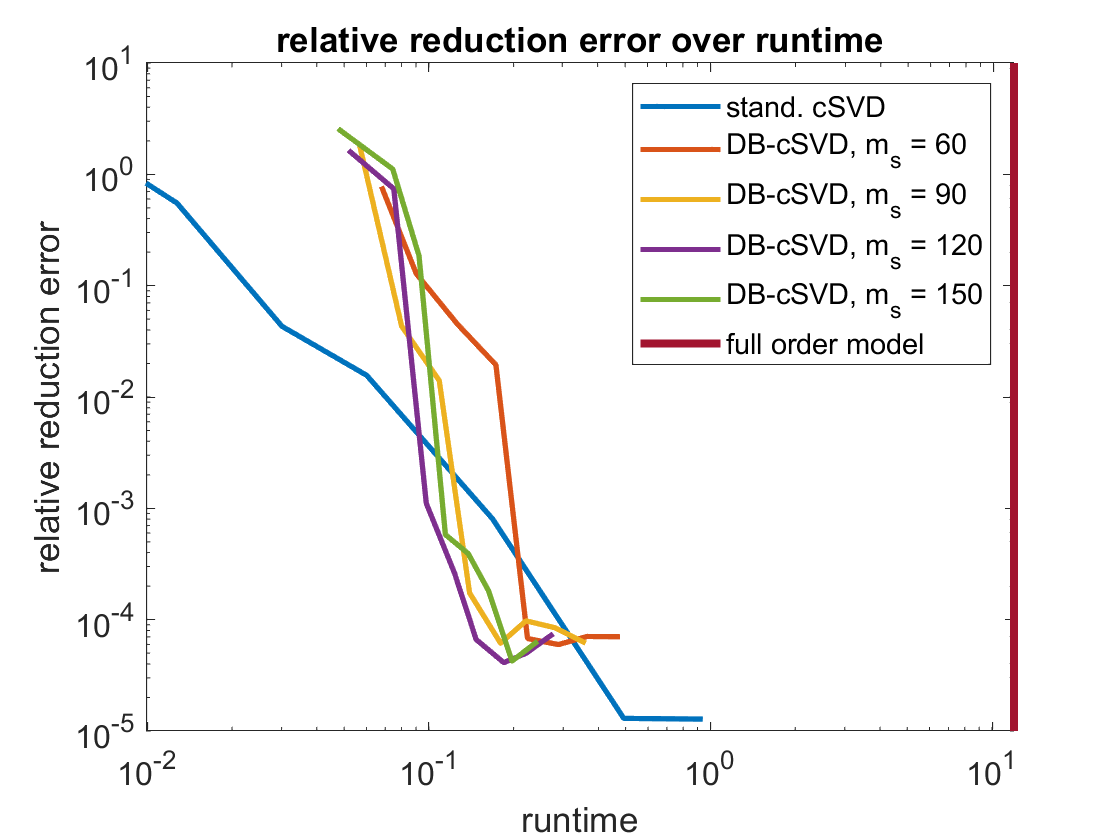}
\end{minipage}
\end{figure}

The results from the reproduction experiments generalize well to unseen data. Repeating the experiments with 10 random parameters from the parameter domain $\mathcal {P}$, window sizes $m_s \in \{60, 90, 120, 150 \}$, numbers of selected snapshots $\ns \in \{100, 150,...,400\}$ and averaging the results for the 10 random parameters leads to the relative reduction errors presented in \Cref{fig:methcomp1gen}. The relative reduction error is again presented in dependence on the average number of basis vectors (see \Cref{avgbasissisze,DBcSVDselcrit}) and in dependence on the online-runtime. We observe again that by using a dictionary-based approach far less basis vectors leads to the same reduction error as a standard approach. For example a relative error of  $10^{-4}$ is obtained with the dictionary-based methods with about 100 basis vectors compared to 500 basis vectors required for a standard cSVD basis. This reduced number of basis vectors also results in better runtimes. A relative error of $10^{-4}$ can be achieved with about half of the runtime compared to the standard approach.

\begin{figure}[H]
\caption{Linear wave equation: Average relative error in Hamiltonian over time-steps, generalization experiment}
\label{fig:methcomp1genHam_non_adapt}
\begin{minipage}[t]{0.475\textwidth}
\includegraphics[width = \textwidth]{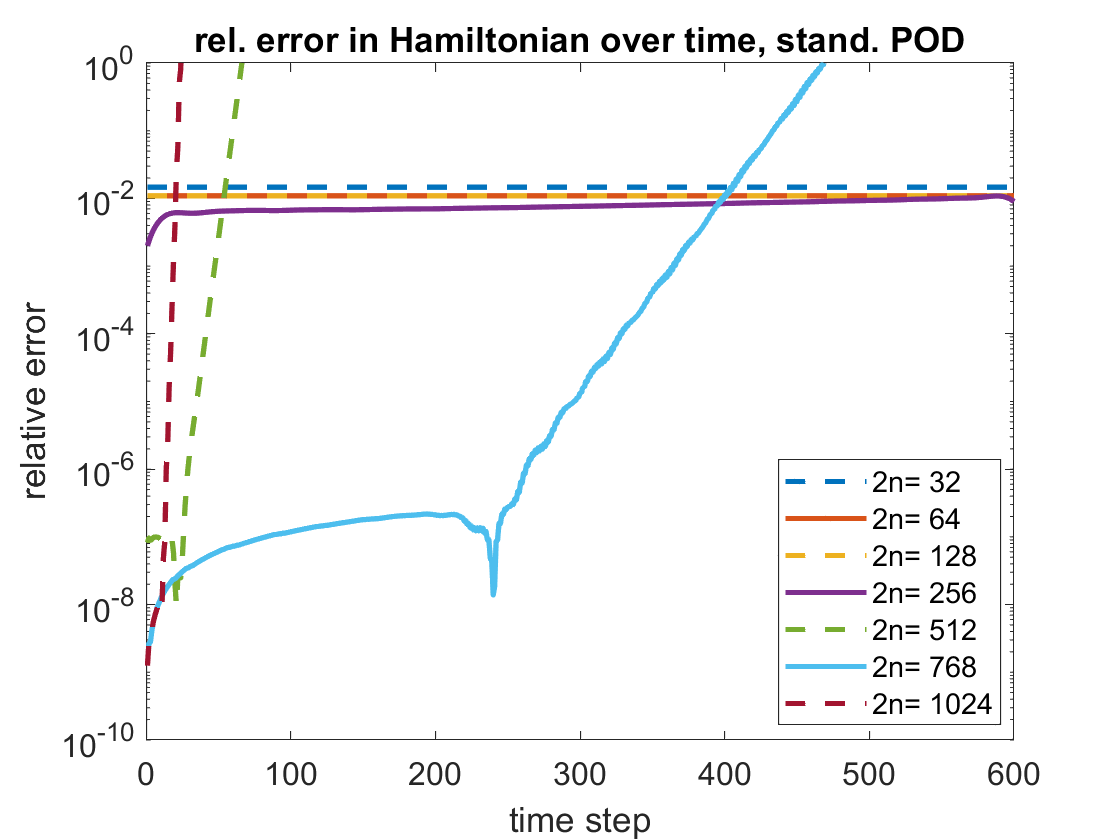}
\end{minipage}
\begin{minipage}[t]{0.475\textwidth}
\includegraphics[width = \textwidth]{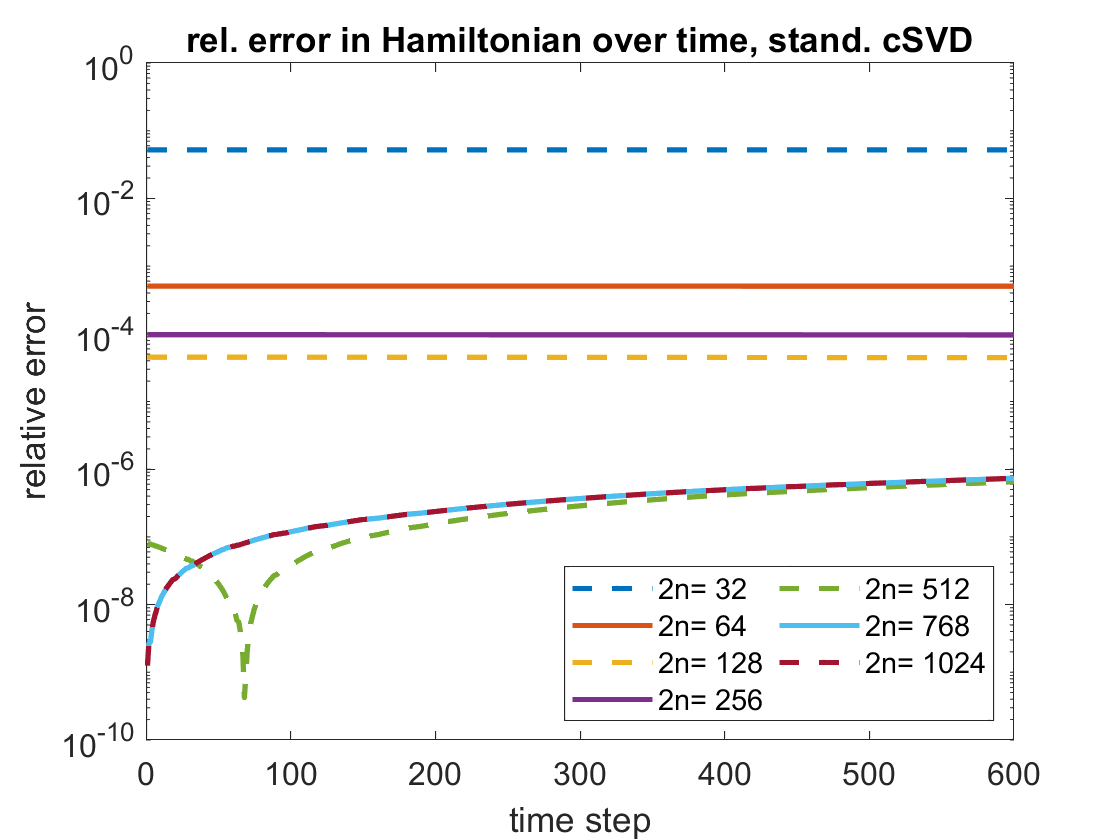}
\end{minipage}
\end{figure}

\begin{figure}[H]
\caption{Linear wave equation: Average relative error in Hamiltonian over time-steps, generalization experiment}
\label{fig:methcomp1genHam}
\begin{minipage}[t]{0.475\textwidth}
\includegraphics[width = \textwidth]{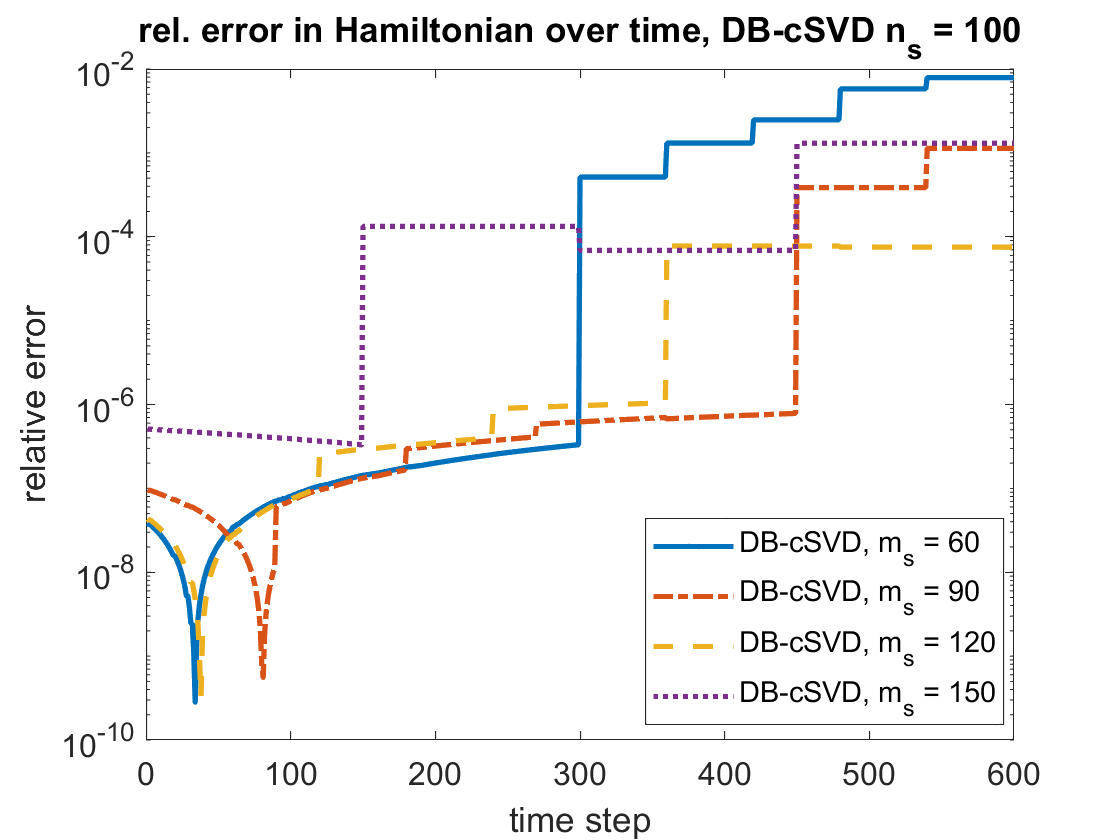}
\end{minipage}
\begin{minipage}[t]{0.475\textwidth}
\includegraphics[width = \textwidth]{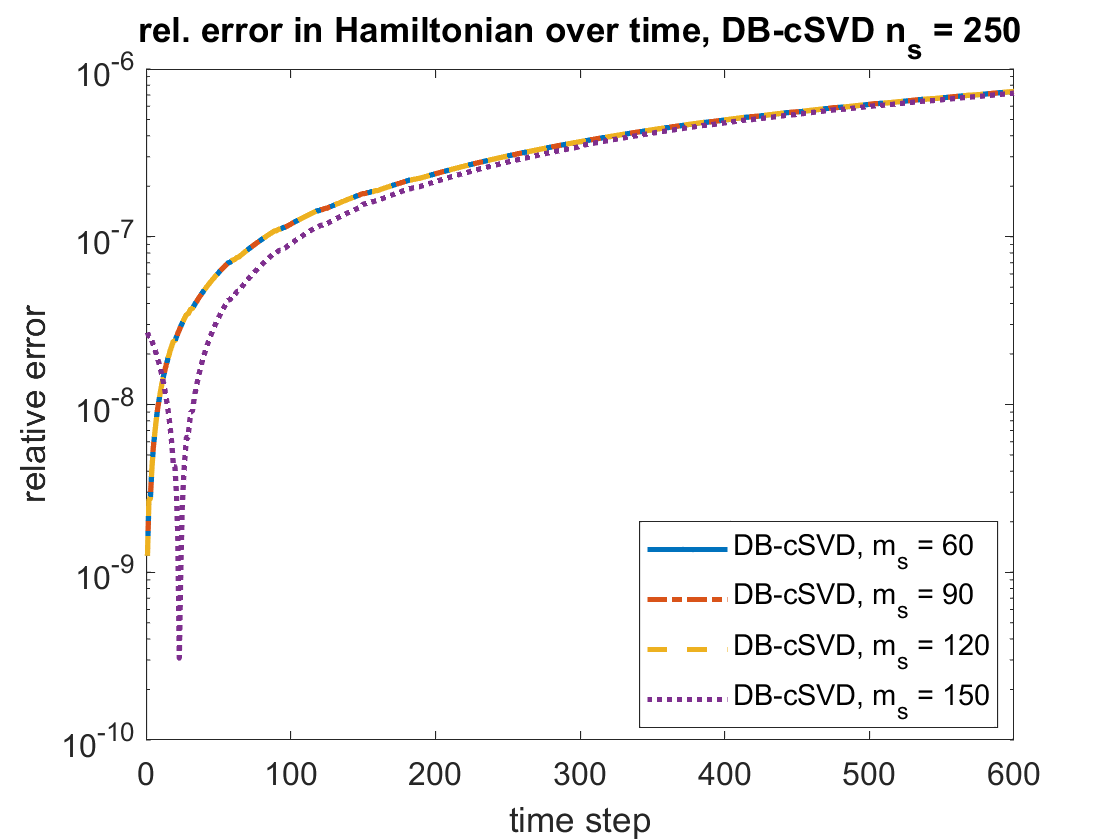}
\end{minipage}
\end{figure}
Also the relative error in the Hamiltonian in dependence on the time-step averaged over 10 random parameters is studied. In \Cref{fig:methcomp1genHam} (analogously to \Cref{fig:methcomp1ham} but for the generalization experiment instead of reproduction), we present the evolution of the error in the Hamiltonian for the standard cSVD. We observe that with 512 or more basis vectors, the relative error in the Hamiltonian is less than $10^{-6}$. Using a non-symplectic, standard POD-basis leads again to high errors in the Hamiltonian for basis sizes $\geq 256$. 
For the DB-cSVD with $n_s=100$ selected snapshots per basis update, we observe high error jumps across the basis changes. We avoid this with $n_s = 250$ selected snapshots and achieve a relative error of less than $10^{-6}$ as in the standard cSVD case. 

\subsection{Sine-Gordon Equation}
\label{sec:sinegordon}
In order to inspect the qualities of the dictionary-based SDEIM, a non-linear wave-equation model is considered, the 1D-Sine-Gordon equation. For the spatial variable  $z\in\varOmega := (0, 50)$ and $t \in  \It(\fmu) =[t_0, \tEnd(\fmu)]$ it reads
\begin{align*}
u_{tt} &=  u_{zz}-s(u) &&\textrm{in } \It(\fmu)\times\varOmega \\ 
u(t_0, z; v) &= u^0(z; v) := 4 \arctan\left(\phi(z;v)\right) &&\textrm{in } \varOmega\\
u_t (t_0, z; v) &= v^0(z; v):= \frac{-4 v}{\sqrt{1-v^2}} \frac{\phi(z;v)} {1 + \phi(z;v)^2} &&\textrm{in } \varOmega\\ 
u(t, 0) & = 0,\ u(t, 50) =2\pi &&\textrm{in } \It(\fmu)
\end{align*}
with
\[\phi(z,v):=\exp\left(\frac{z - 10} {\sqrt{1-v^2}}\right)\ \textrm{and}\ s(u) = \sin (u).\]
We choose $\tInit = 0,\, \tEnd(\fmu) = 30/\fmu$ and as parameter (vector) $\fmu = v\in \mathcal{P} := [0.7,0.9] $. 
In order to obtain homogeneous Dirichlet boundary conditions we reformulate the system. This has the advantage that the boundary condition will exactly be conserved by a reduced simulation because then all state snapshots will have zero displacement on the boundary. This inherits to the reduced solution, which is an element of the span of all state snapshots. 

We split the solution $u$ in a function $u$ and the initial condition $u^0$. 
\[u(t,z)= \hat u(t,z) + u^0(z).\] 
Inserting this into the PDE 
\[u_{tt} = c^2 u_{zz} - s(u)\]
leads to
 \[ \hat u_{tt} + u^0_{tt} = c^2 ( \hat u_{zz} + u^0_{zz}) - s( \hat u + u^0).\]
As $ (u^0)_{tt} = 0$ and 
\[u^0_{zz}  = \frac{2}{v^2-1}\frac{\textrm{tanh}(\frac{z-10}{\sqrt{1-v^2}})}{\textrm{cosh}(\frac{z-10}{\sqrt{1-v^2}})},\]
we get
\[\hat u_{tt}  = c^2  \hat u_{zz} -\frac{2c^2}{1-v^2}\frac{\textrm{tanh}\left(\frac{z-10}{\sqrt{1-v^2}}\right)}{\textrm{cosh}\left(\frac{z-10}{\sqrt{1-v^2}}\right)} - s( \hat u +  u^0).\]

For the 1D-Sine-Gordon equation, a spatial central finite difference discretization of second order with the positive definite matrix $\fD_{zz}\in \R^{N \times N}$ and grid points $z_1..., z_N$ leads to
\begin{align*}\label{sine_discr}
\ddt \fx(t; \fmu)  &= \JtN \grad[\fx]\Ham(\fx(t; \fmu)) = \JtN\fH\fx(t;\fmu)+ \JtN\ff_{\textrm{nl}}(\fx(t; \fmu)) + \JtN \fb,\\ 
\fx(0; \fmu) &=  \fxInit(\fmu)
\end{align*}
with 
$$\fH= \begin{pmatrix}
\fD_{zz} \ &\Z{N} \\
\Z{N}\  &\I{N}
\end{pmatrix}, \ff_{\textrm{nl}}(\fx) =  \bigg{[}s(x_1+u^0(z_1));...; s(x_N+u^0(z_N)); \Z{N\times 1}\bigg{]},$$
parametric initial value
$$\fxInit(\fmu) = [\Z{N\times 1};v^0(z_1; \fmu)); ...,v^0(z_N; \fmu))] $$
and
\[\fb = \left[\frac{2}{v^2-1}\frac{\textrm{tanh}\left(\frac{z_1-10}{\sqrt{1-v^2}}\right)}{\textrm{cosh}\left(\frac{z_1-10}{\sqrt{1-v^2}}\right)};...;\frac{2}{v^2-1}\frac{\textrm{tanh}\left(\frac{z_N-10}{\sqrt{1-v^2}}\right)}{\textrm{cosh}\left(\frac{z_N-10}{\sqrt{1-v^2}}\right)}; \Z{N\times 1}\right] \in \R^{2N}. \]\\  \\ 
The corresponding Hamiltonian reads 
\[\Ham(\fx) = \frac{1}{2}\rT\fx\fH\fx  + \rT{\mathds{1}_{N}}\overline \fS(\fx)  + \rT\fb\fx\] 
with 
$$\overline \fS (\fx) =   [S(x_1+u^0(z_1));...; S(x_N+u^0(z_N))]\in \R^{N}$$
for some function $S$ with   $S'(u) = s(u).$ The number of grid points is chosen as $N_{z}= 5000$, which leads to a size of $2N = 10000$ of the Hamiltonian system. The number of equidistant time-steps is chosen as $n_t = 400$. Again, this results in different time-step sizes for different parameters. First, we study a reproduction experiment. For the dictionaries, we use the snapshots from the parameters $\fmu_j = 0.7 + (j-1)\cdot0.05, j = 1,...,5$ from which also the snapshot matrix for the standard cSVD-basis is computed. This results in a dictionary size of $\NX = 2000.$
Now that we know the structure of $\ff_{\textrm{nl}}$, we explain how the selection $\rT\fP\ff_{\textrm{nl}}(\fV\fy)$ is calculated online-efficiently.
Because each component of 
$$\ff_{\textrm{nl}}(\fx)=  \bigg{(}s(x_1+u^0(z_1));...; s(x_N+u^0(z_N)), \Z{1\times N}\rT{\bigg{)}}$$
depends only either on one entry or no entry of $\fx$, the DEIM-algorithm is efficiently applicable. It follows, that the function evaluation 
$\ff_{\textrm{nl}}(\fV\fy)$ depends only on a few entries of $\fx$. 
The $i$th entry of $\rT\fP\ff_{\textrm{nl}}(\fV\fy)$ is
\begin{equation}\label{selectinds}
(\rT\fP\ff_{\textrm{nl}}(\fV\fy))_i =
\begin{cases}
  s\big((\fV\fy)_{\rho_i}\big),
& \rho_i \leq N,\\
  0, &  \rho_i > N,
\end{cases}
\end{equation}
with $\rho_i$ the DEIM indices from \DEIMalg. Therefore, $(\rT\fP\ff_{\textrm{nl}}(\fV\fy))$ can be computed as
\begin{equation}\label{DEIMOnlEff}
(\rT\fP\ff_{\textrm{nl}}(\fV\fy)) = (\ff_{\textrm{nl}}(\rT\fP\fV\fy))
\end{equation}
because no $\rho_i > N$ will be chosen by the \DEIMalg-algorithm as $\ff_{\textrm{nl}}(\fx)_i = 0$ for $i > N.$
During DB-SDEIM, the projection $\rT\fP_{\textrm{o}}\rT{\hat\fP_{\dictP}}\ff_{\textrm{nl}}(\fV\fy)$ is computed in two steps. In the offline-phase $\hat\fP_{\dictP}$ is precomputed and during the reduced simulation $[\frho_{\textrm{o}},\fP_{\textrm{o}}]=$ \DEIMalg($\fFhatPs\widetilde \fPsi$) is computed. According to \Cref{DEIMOnlEff}, the projection $\rT\fP_{\textrm{o}}\rT{\hat\fP_{\dictP}}\ff_{\textrm{nl}}(\fV\fy)$ can be computed as $\rT\fP_{\textrm{o}}\rT{\hat\fP_{\dictP}}\ff_{\textrm{nl}}(\fV\fy)$  = $\ff_{\textrm{nl}}(\rT\fP_{\textrm{o}}\rT{\hat\fP_{\dictP}}\fV\fy)$. Therefore, the goal is to assemble the product $\rT\fP_{\textrm{o}}\rT{\hat\fP_{\dictP}}\fV\fy$ without explicitly computing the matrices $\fP_{\textrm{o}}, \hat\fP_{\dictP}$ and $\fV$. For that we use the representation $\fV =[\fY[\fP_s; \fP_s] \widetilde \fPhi ,\TJtN \fY[\fP_s; \fP_s] \widetilde \fPhi]$ from \Cref{tablestandvsdictbasedcSVD}, with $\fY= [\fX,\JtN\fX]$. This yields $$\fV =[[\fX\fP_s,\JtN\fX\fP_s]\widetilde \fPhi ,[-\JtN\fX\fP_s, \fX\fP_s]\widetilde \fPhi]$$  using the facts, that $\TJtN = -\JtN$ and $\JtN\TJtN=\TJtN\JtN = \fI_{2N}$. 
In the offline-phase we precompute 
$$\fG_{\hat\fP_{\dictP},\textrm{X}}= \rT{\hat\fP_{\dictP}}\fX\ \ \textrm{and}\ \ \fG_{\hat\fP_{\dictP},\JtN,\textrm{X}}= \rT{\hat\fP_{\dictP}}\JtN\fX$$ via index selection. Then, during the online-phase 
$$\fG_{\hat\fP_{\dictP},\textrm{X}, \textrm{s}} =\fG_{\hat\fP_{\dictP},\textrm{X}}\fP_s\ \ \textrm{and}\ \ \fG_{\hat\fP_{\dictP},\JtN,\textrm{X}, \textrm{s}} =\fG_{\hat\fP_{\dictP},\JtN,\fXs}\fP_s$$ are calculated via index selection. Finally, the projections 
$$\fG_{\fP_{\textrm{o}},\hat\fP_{\dictP},\textrm{X}, \textrm{s}} = \rT\fP_{\textrm{o}}\fG_{\hat\fP_{\dictP},\textrm{X}, \textrm{s}}\ \ \textrm{and}\ \ \fG_{\fP_{\textrm{o}},\hat\fP_{\dictP},\JtN,\textrm{X}, \textrm{s}} = \rT\fP_{\textrm{o}} \fG_{\hat\fP_{\dictP},\JtN,\textrm{X}, \textrm{s}}$$ 
are computed and from $\fG_{\hat\fP_{\dictP},\JtN,\textrm{X}, \textrm{s}}\widetilde \fPhi$ and $\fG_{\fP_{\textrm{o}},\hat\fP_{\dictP},\JtN,\textrm{X}, \textrm{s}}\widetilde \fPhi$ the product 
$$\rT\fP_{\textrm{o}}\rT{\hat\fP_{\dictP}}\fV = \rT\fP_{\textrm{o}}\rT{\hat\fP_{\dictP}}[[\fX\fP_s,\JtN\fX\fP_s]\widetilde \fPhi ,[\TJtN\fX\fP_s, \fX\fP_s]\widetilde \fPhi]$$
is stacked.
The size of the $i$th DB-SDEIM-basis $\widetilde m_i$ is chosen accordingly to the size of the DB-cSVD-basis (see \Cref{DBcSVDselcrit}) as the lowest number $\widetilde m$ that fulfills
\begin{equation}\label{DBSDEIMselcrit}
(1-\epsilon_\textrm{SDEIM})\sum\limits_{\ell = 0}^{n_s}\lambda_{\ell}^i<\sum\limits_{\ell = 0}^{\widetilde m+1}\lambda_{\ell}^i,
\end{equation}
where we denote with $\lambda_{\ell}^i>0, \ell = 1,..., \widetilde m$ the dominant non-zero eigenvalues of $\fGFs^i := (\fP_s^i\rT) \fGF\fP_s^i$, with $\fP_s^i$ the $i$th selection matrix. For the Sine-Gordon experiments we choose $\epsilon_\textrm{SDEIM}=10^{-12}$ and $\epsilon_\textrm{cSVD}=10^{-13}.$

\begin{figure}[H]
\caption{Non-linear Sine-Gordon: Average relative reduction error over number of basis vectors and online-runtime, reproduction experiment}
\label{fig:methcomp2rep}
\begin{minipage}[t]{0.475\textwidth}
\includegraphics[width = \textwidth]{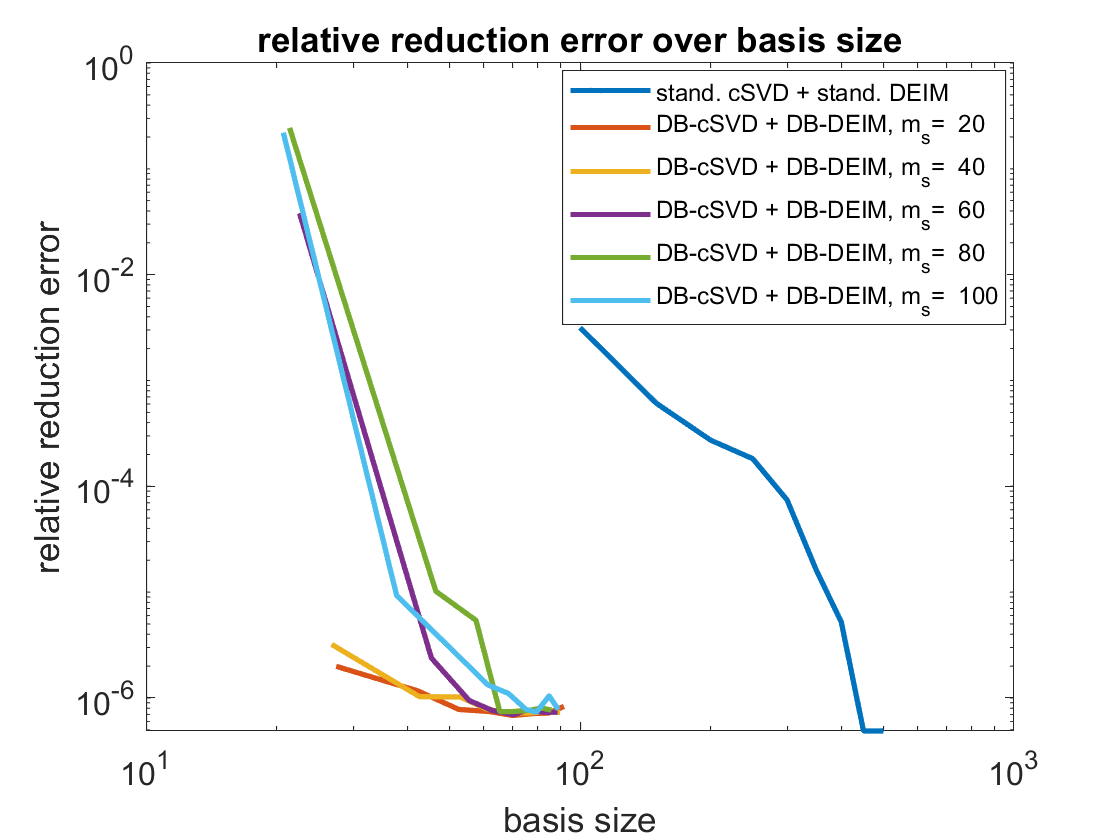}
\end{minipage}
\begin{minipage}[t]{0.475\textwidth}
\includegraphics[width = \textwidth]{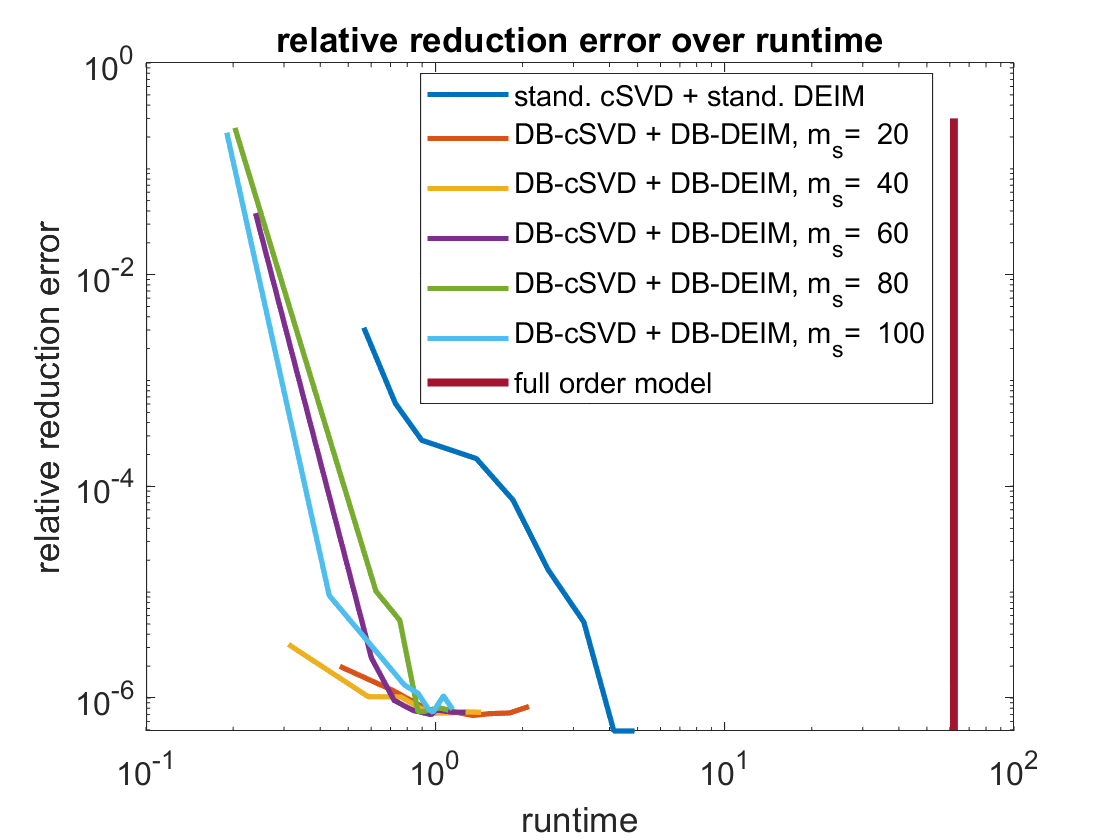}
\end{minipage}
\end{figure}

In \Cref{fig:methcomp2rep} we present the relative reduction error $\fe_\text{rel}(\fmu)$ (from \Cref{eqnrelerr}) in dependence on the average basis size and in dependence on the online-runtime. All reduction errors are averaged over the 5 training parameters. For the DB-cSVD on the horizontal axis, the average number of the basis vectors (see \Cref{avgbasissisze,DBcSVDselcrit}) is presented. The DB-cSVD with DB-SDEIM is compared for different window sizes $m_s \in \{20, 40, 60, 80, 100\}$ and numbers of selected snapshots $\ns \in \{50, 100, 150,...,400\}$ to the standard cSVD with standard SDEIM and a standard POD with standard DEIM. 
Similarly to the experiments for the linear wave equation, using a dictionary-based approach, the same reduction error is achieved by smaller average basis sizes. For example a relative error of  $10^{-6}$ is obtained using in average $40-50$ basis vectors with the dictionary-based methods compared to 400 basis vectors required for a standard cSVD basis. This reduced number of basis vectors also leads to better runtimes. With the dictionary-based methods a relative error of $10^{-6}$ is achieved with an online-runtime of less than one second compared to about 4 seconds with the standard approach. 

In \Cref{fig:methcomp2repham,{fig:methcomp2rephamadapt}} we present the relative error in the Hamiltonian 
$\fe_{\Ham, \text{rel}, i}(\fmu)$ (see Equation \Cref{eqn:relHam}) averaged over the five training parameters in dependence on the time-step.
\begin{figure}[H]
\caption{Non-linear Sine-Gordon: Average relative error in Hamiltonian over time-steps, reproduction experiment}
\label{fig:methcomp2repham}
\begin{minipage}[t]{0.475\textwidth}
\includegraphics[width = \textwidth]{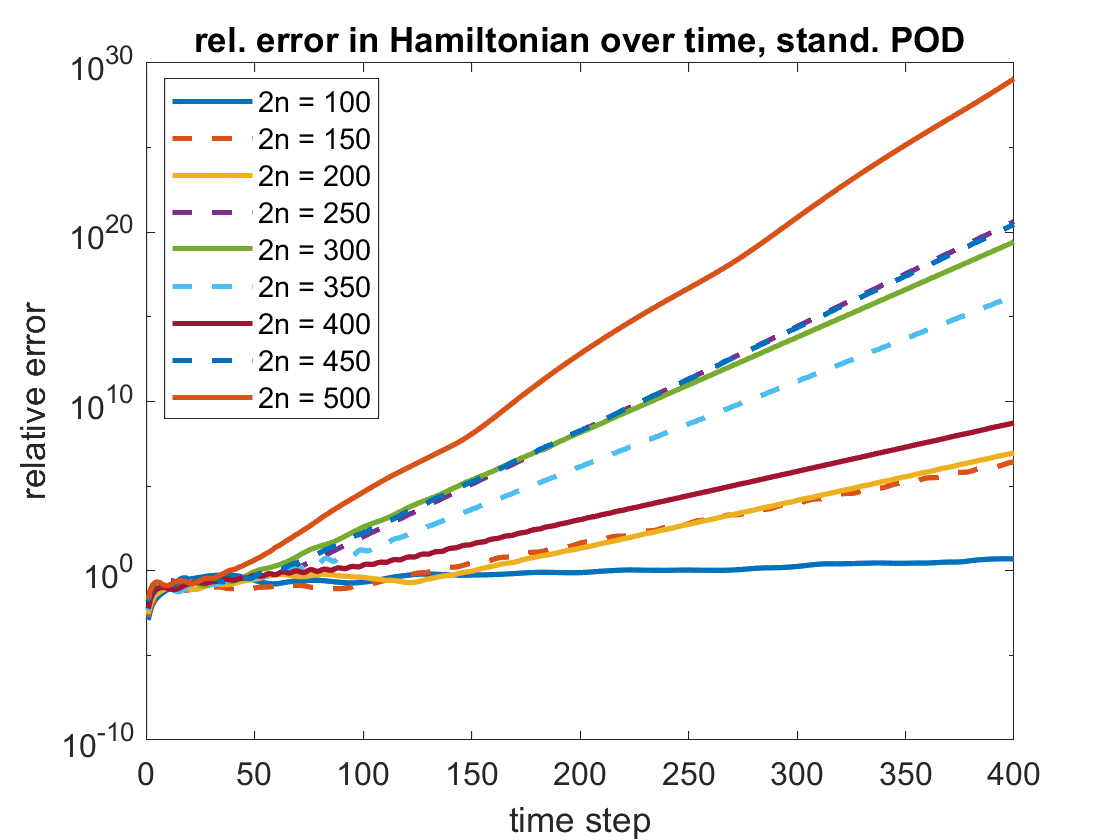}
\end{minipage}
\begin{minipage}[t]{0.475\textwidth}
\includegraphics[width = \textwidth]{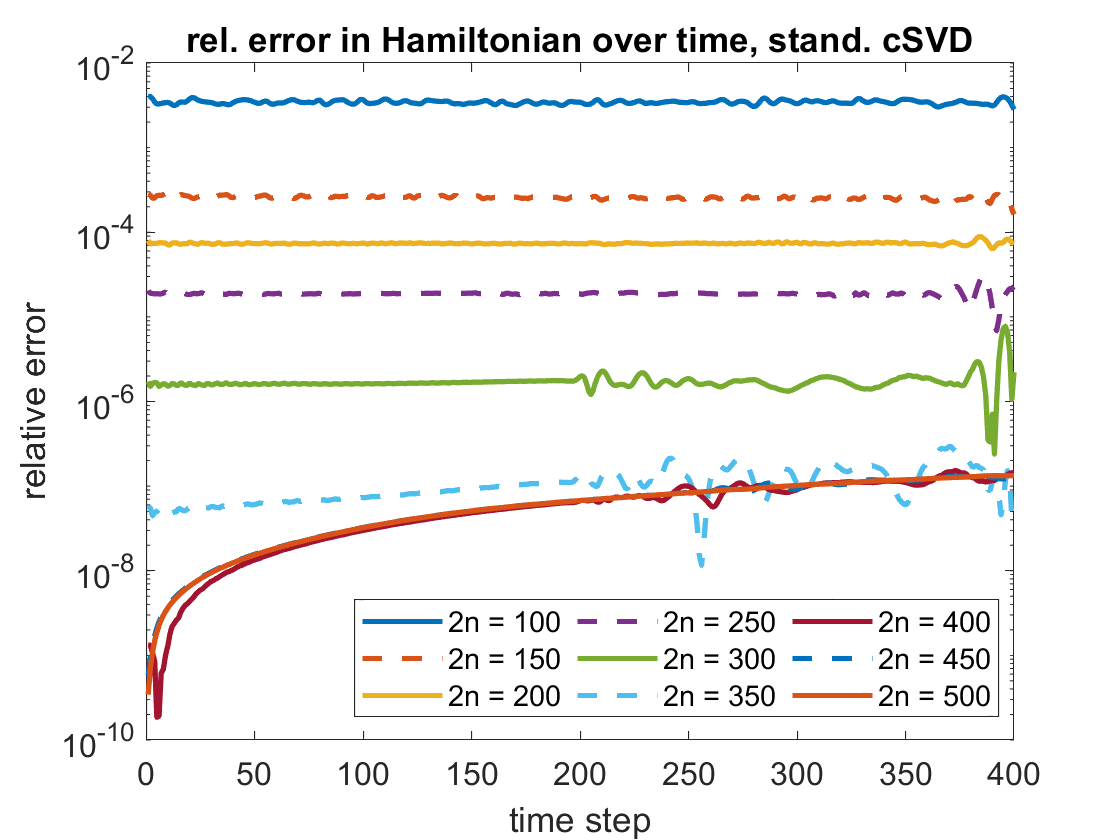}
\end{minipage}
\end{figure}

\begin{figure}[H]
\caption{Non-linear Sine-Gordon: Average relative error in Hamiltonian over time-steps, reproduction experiment}
\label{fig:methcomp2rephamadapt}
\begin{minipage}[t]{0.475\textwidth}
\includegraphics[width = \textwidth]{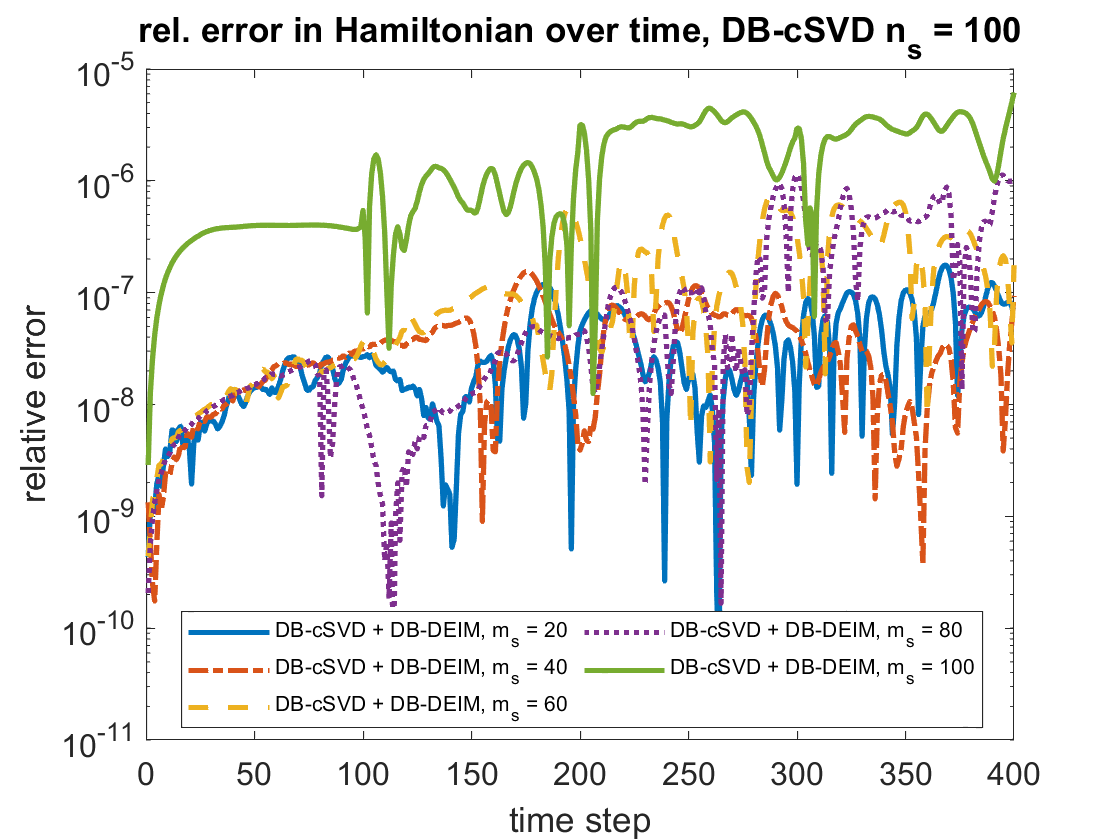}
\end{minipage}
\begin{minipage}[t]{0.475\textwidth}
\includegraphics[width = \textwidth]{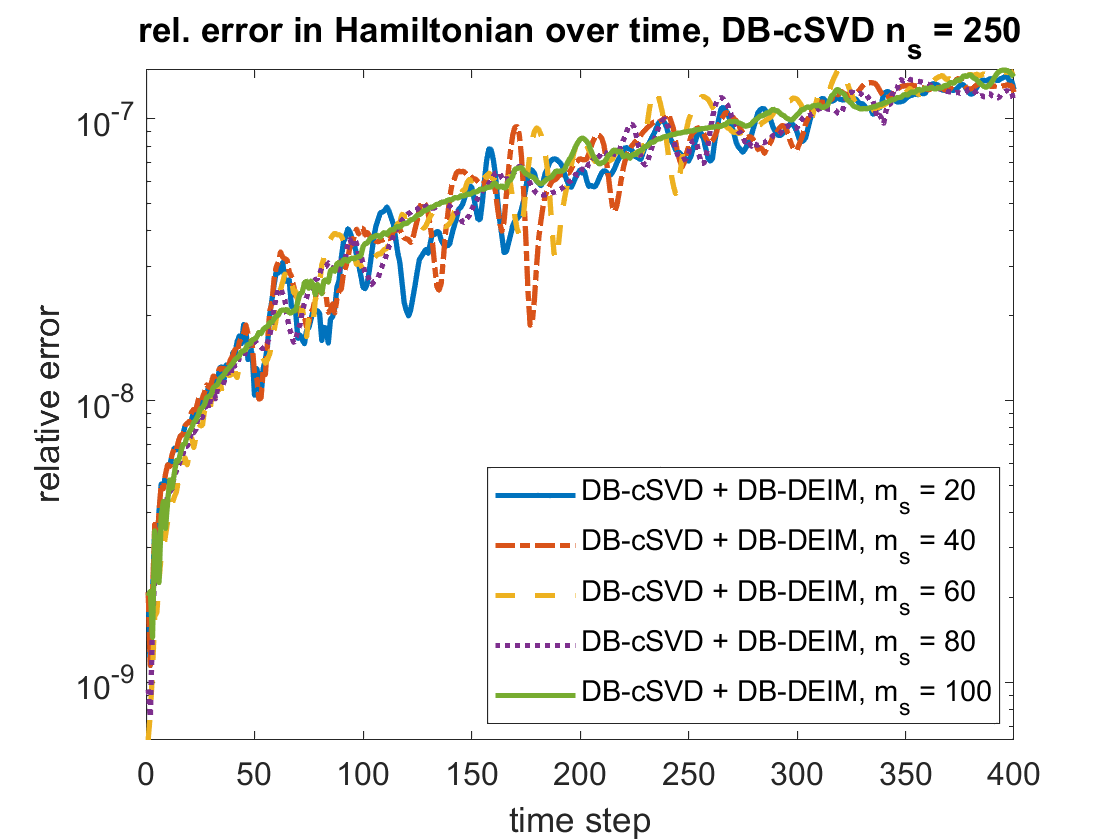}
\end{minipage}
\end{figure}
We observe that with 350 or more basis vectors, the relative error in the Hamiltonian is smaller than $10^{-8}$ for the standard approach. Note that we do not expect exact preservation of the reduced Hamiltonian and thus non-constant error in the Hamiltonian as the time stepping only preserves Hamiltonian functions up to quadratic order. 
With the dictionary-based approach and $\ns = 250$ we achieve a relative error of about $10^{-7}$. 

\begin{figure}[H]
\caption{Non-linear Sine-Gordon: Average relative reduction error over number of basis vectors and online-runtime, generalization experiment}
\label{fig:methcomp2gen}
\begin{minipage}[t]{0.475\textwidth}
\includegraphics[width = \textwidth]{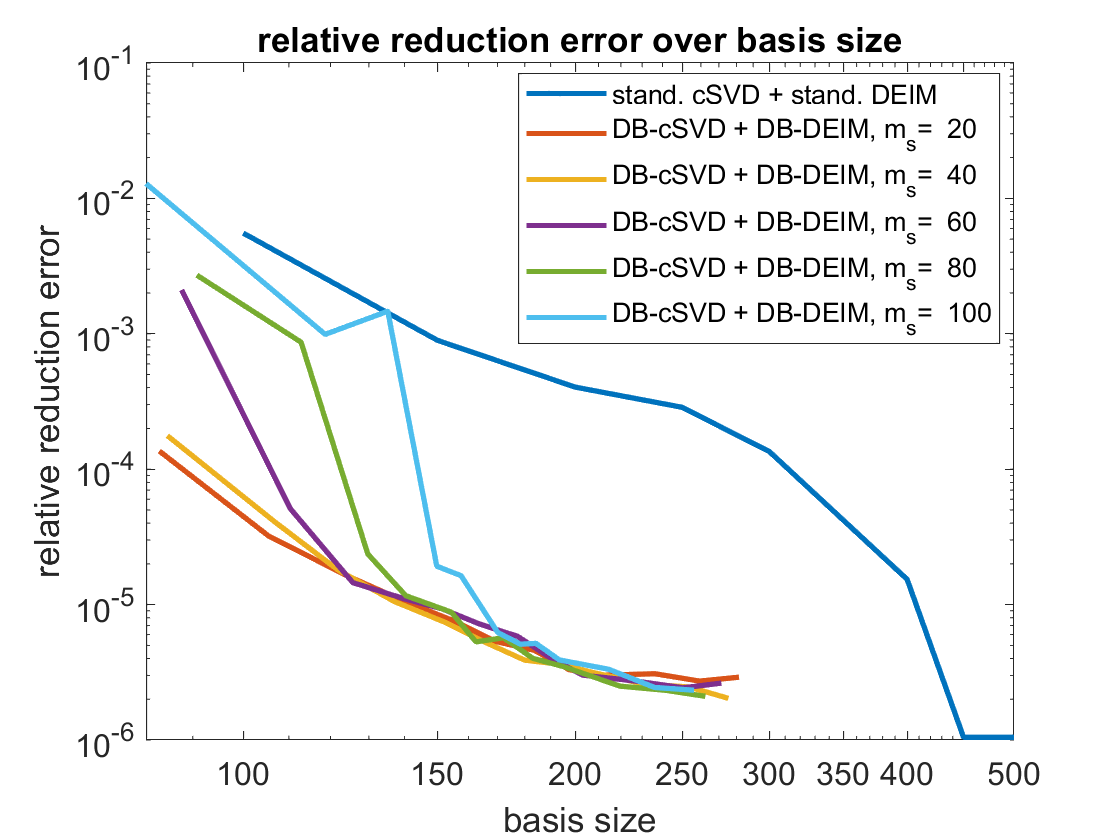}
\end{minipage}
\begin{minipage}[t]{0.475\textwidth}
\includegraphics[width = \textwidth]{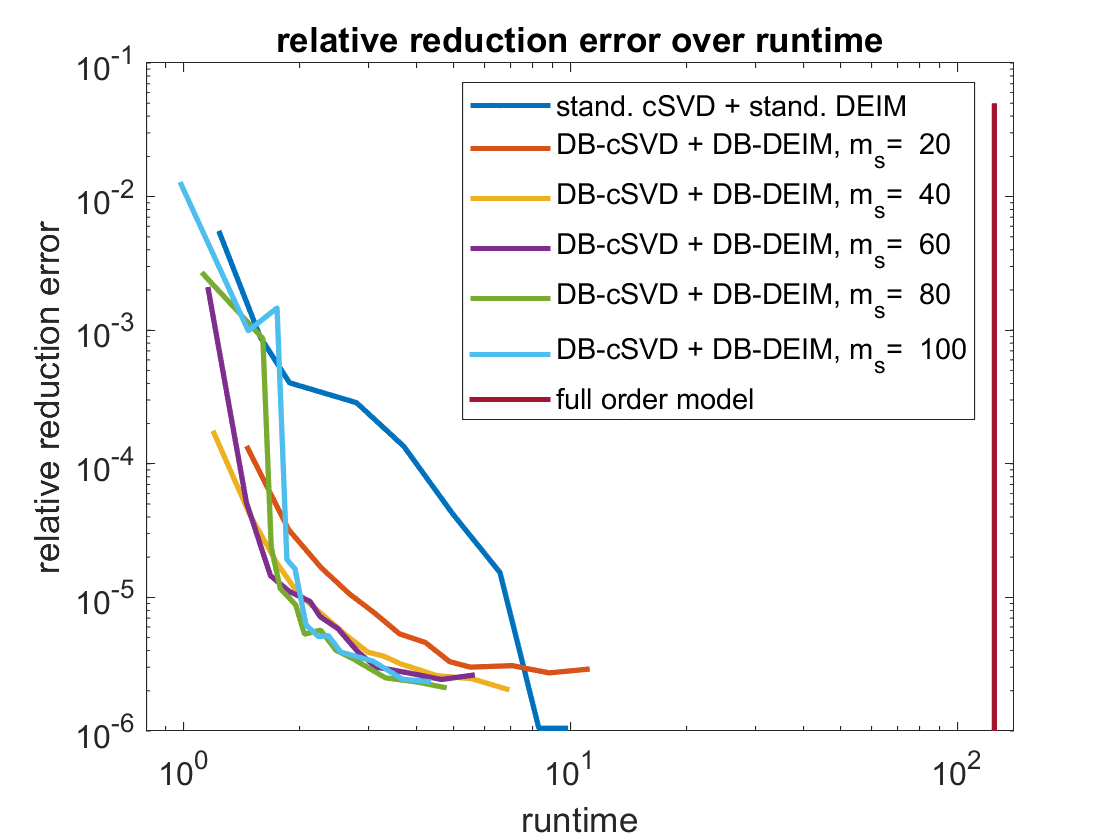}
\end{minipage}
\end{figure}

The results from the reproduction experiments generalize well  to unseen data. We repeat the experiments with 10 random parameters, different window sizes $m_s \in \{20, 40, 60, 80, 100\}$ and numbers of selected snapshots $\ns \in \{100, 150,...,500\}$ and average the results. This leads to the relative reduction errors presented in \Cref{fig:methcomp2gen}. The relative reduction error is again presented in dependence on the average number of basis vectors (see \Cref{avgbasissisze,DBcSVDselcrit}) and in dependence on the online-runtime.

With the dictionary-based methods a relative error of  for example $10^{-5}$ is achieved with about 140 basis vectors. For a standard cSVD basis 500 basis vectors have been required. This reduced number of basis vectors also results in better runtimes. With the dictionary-based approach we obtain a relative error of $10^{-5}$ with an online-runtime of less than two seconds compared to about 8 seconds required for the same error with the standard approach. 

\begin{figure}[H]
\caption{Non-linear Sine-Gordon: Average relative error in Hamiltonian over time-steps, generalization experiment}
\label{fig:methcomp2genhamnonadapt}
\begin{minipage}[t]{0.475\textwidth}
\includegraphics[width = \textwidth]{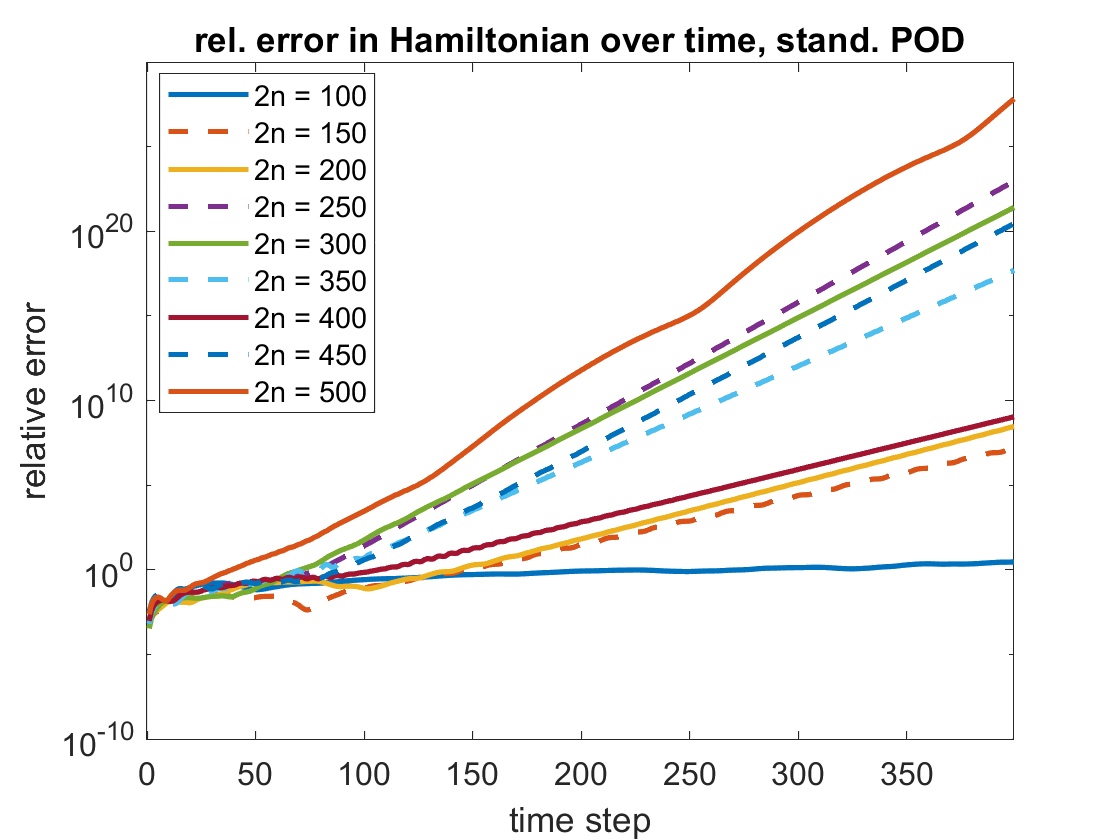}
\end{minipage}
\begin{minipage}[t]{0.475\textwidth}
\includegraphics[width = \textwidth]{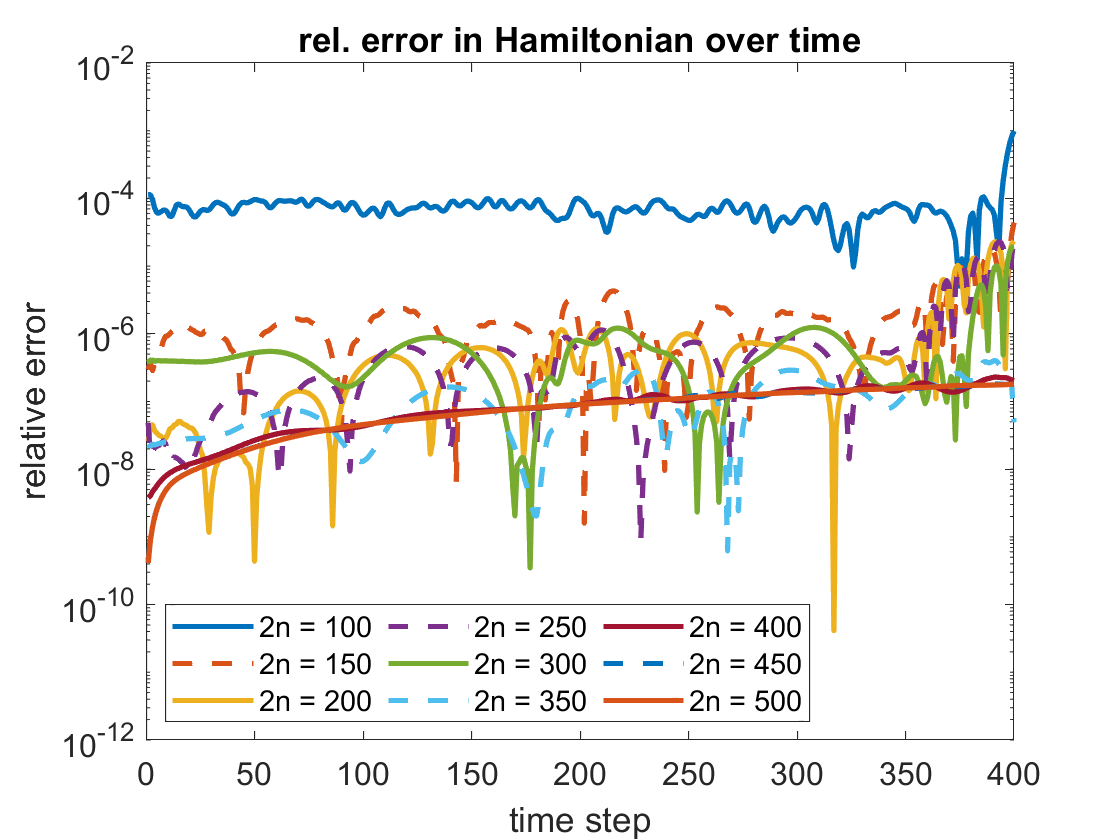}
\end{minipage}
\end{figure}

\begin{figure}[H]
\caption{Non-linear Sine-Gordon: Average relative error in Hamiltonian over time-steps, generalization experiment}
\label{fig:methcomp2genham}
\begin{minipage}[t]{0.475\textwidth}
\includegraphics[width = \textwidth]{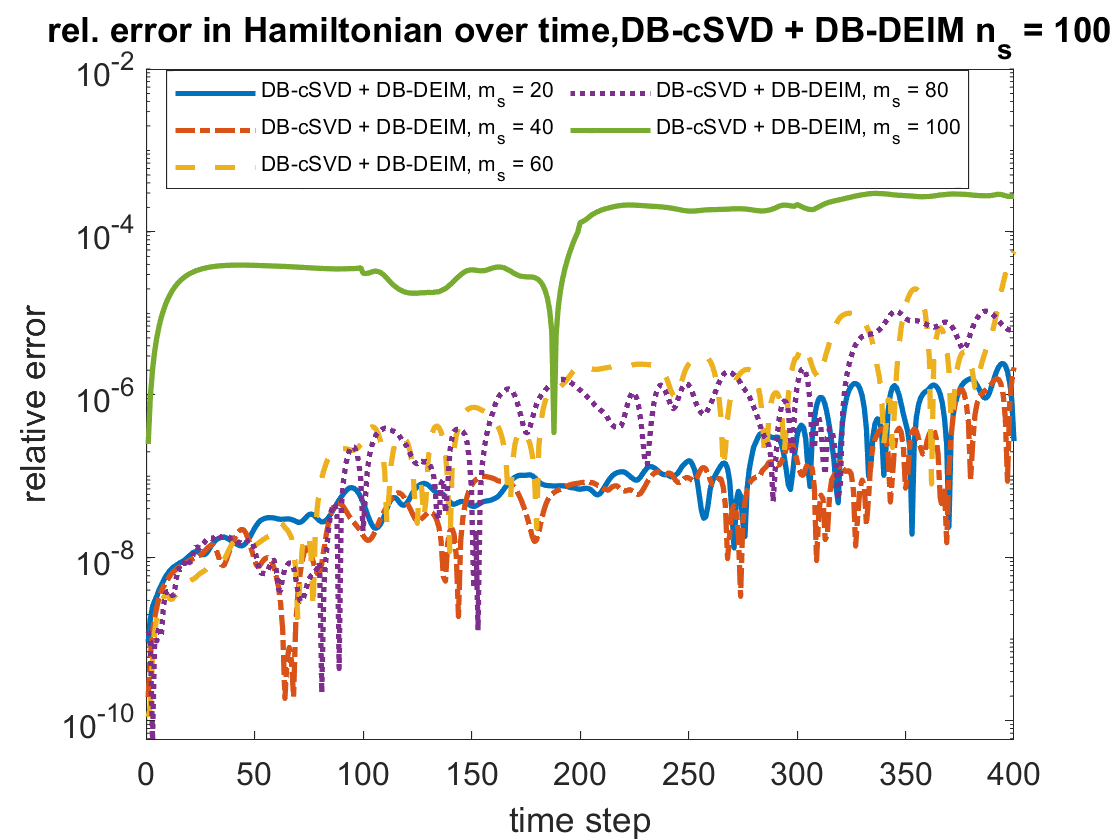}
\end{minipage}
\begin{minipage}[t]{0.475\textwidth}
\includegraphics[width = \textwidth]{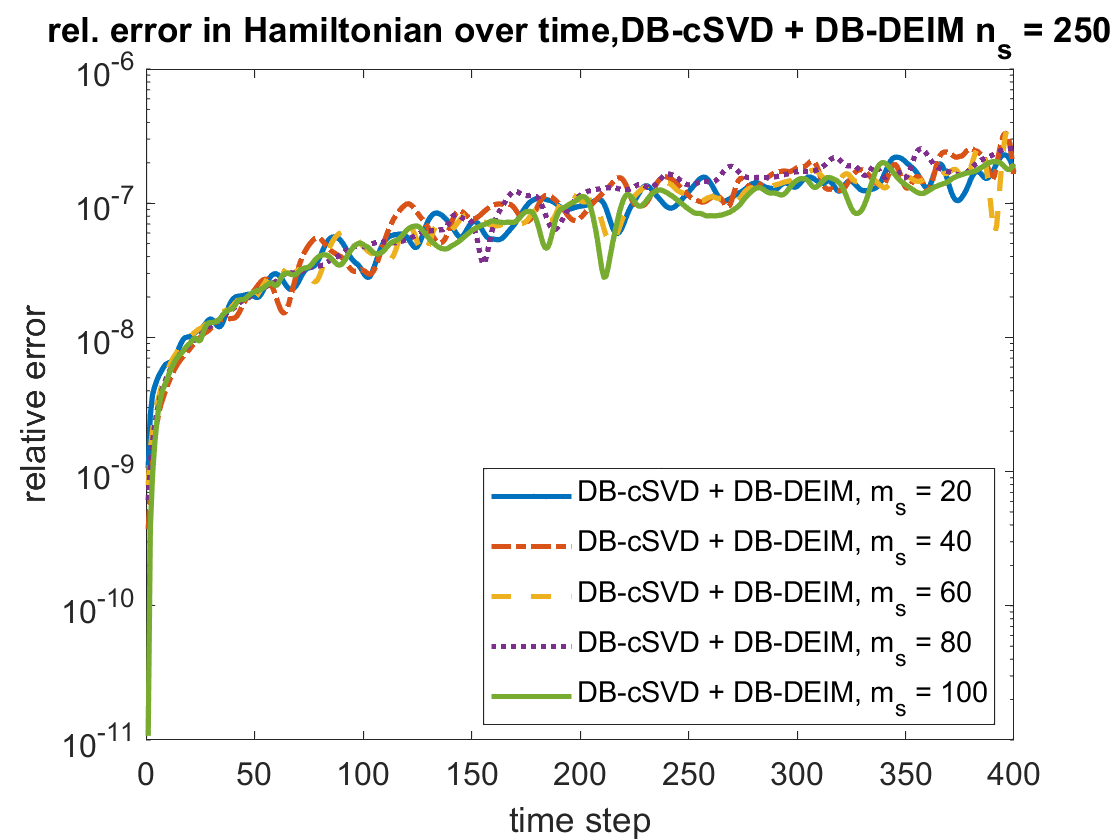}
\end{minipage}
\end{figure}

In \Cref{fig:methcomp2genhamnonadapt,fig:methcomp2genham} we present the average relative error in the Hamiltonian in dependence on the time-step.

We observe that with the dictionary-based approaches (nearly) the same relative error in the Hamiltonian is achieved. 
With 350 or more basis vectors the relative error in the Hamiltonian is about $10^{-7}$ for the standard cSVD with standard SDEIM and for the dictionary-based version with $\ns = 250$.

\section{Conclusion and Outlook}\label{sec:conclusion}
In this work we presented a new dictionary-based, symplectic MOR technique, the DB-cSVD. Furthermore, we introduced new dictionary-based hyper-reduction techniques for non-linear systems, the DB-DEIM and DB-SDEIM. We derived offline-online-splitting techniques that allow us to perform the computations in the online-phase to be independent of the state dimension and to be only dependent on the size of the dictionary and the number of selected snapshots $\ns.$ Furthermore, we derived an error bound for the error in the Hamiltonian resulting from the projections during the reduced simulation. In the numerical experiment section, we showed, that in practice the error in the Hamiltonian computed with the dictionary-based methods is at most slightly higher than the error obtained by standard symplectic methods if the number of snapshots is chosen high enough. For the relative reduction error we showed, that the same error can be obtained by far smaller (average) basis sizes, which also resulted in corresponding speed-ups. 

There exist several open questions and options for future work. For example, 
the dictionary is so far constructed by fixed a priori sampling of the parameter-time domain. As the dictionary size directly enters into the computational complexity of the online-scheme, a more adaptive generation of the dictionary would presumably enable a more compact representation. Currently, the selection process is realized by parameter-time-distances. It would be interesting to devise residual-based selection criteria. The difficult point here is the computational efficiency. Furthermore, so far computing the parameter-time-distances a 2-norm is used, where the time-distance is scaled. Generalizing more problem-dependent anisotropic distances as presented in \cite{stam13} to time-dependent problems could presumably lead to an online-speed-up of the dictionary-based methods, as the necessary number of selected snapshots may be reduced.

\section*{Statements and Declarations}

\section*{Funding} 
Funded by Deutsche Forschungsgemeinschaft (DFG, German Research Foundation)
under Project No.~314733389 and Germany’s Excellence Strategy - EXC 2075 – 390740016.
We acknowledge support by the Stuttgart Center for Simulation Science (SimTech).
\section*{Competing interests}  The authors declare that they have no competing interests.

\bibliographystyle{siam}
\bibliography{references_sympl_dicionary_mor}

\end{document}